\documentclass[12pt,a4paper,titlepage]{article}
\usepackage{t1enc}
\usepackage{amsmath}
\usepackage{amssymb}
\usepackage{latexsym}
\usepackage{natbib}
\usepackage{epsfig}
\usepackage{textcomp}
\usepackage{alltt}
\usepackage{graphicx}
\usepackage{rotating}
\usepackage{dcolumn}
\usepackage{xypic,t1enc}
\usepackage{float}
\usepackage{alltt}
\usepackage{color}

\begin{document}
\newcommand{\cip}{\perp\!\!\!\perp}
\newcommand{\nothere}[1]{}
\newcommand{\noi}{\noindent}
\newcommand{\mbf}[1]{\mbox{\boldmath $#1$}}
\newcommand{\cond}{\, |\,}
\newcommand{\hO}[2]{{\cal O}_{#1}^{#2}}
\newcommand{\hF}[2]{{\cal F}_{#1}^{#2}}
\newcommand{\tl}[1]{\tilde{\lambda}_{#1}^T}
\newcommand{\la}[2]{\lambda_{#1}^T(Z^{#2})}
\newcommand{\I}[1]{1_{(#1)}}
\newcommand{\cd}{\mbox{$\stackrel{\mbox{\tiny{\cal D}}}{\rightarrow}$}}
\newcommand{\cp}{\mbox{$\stackrel{\mbox{\tiny{p}}}{\rightarrow}$}}
\newcommand{\cas}{\mbox{$\stackrel{\mbox{\tiny{a.s.}}}{\rightarrow}$}}
\newcommand{\ld}{\mbox{$\; \stackrel{\mbox{\tiny{def}}}{=3D} \; $}}
\newcommand{\nk}{\mbox{$n \rightarrow \infty$}}
\newcommand{\con}{\mbox{$\rightarrow $}}
\newcommand{\dprime}{\mbox{$\prime \vspace{-1 mm} \prime$}}
\newcommand{\Borel}{\mbox{${\cal B}$}}
\newcommand{\bevis}{\mbox{$\underline{\em{Proof}}$}}
\newcommand{\Rd}[1]{\mbox{${\Re^{#1}}$}}
\newcommand{\il}[1]{{\int_{0}^{#1}}}
\newcommand{\pl}[1]{\mbox{\bf {\LARGE #1}}}
\newcommand{\real}{\mathbb{R}}
\newcommand{\bth}{\mbox{\boldmath{$\theta$}}}
\newcommand{\bphi}{\mbox{\boldmath{$\phi$}}}
\newcommand{\bbe}{\mbox{\boldmath{$\beta$}}}
\newcommand{\bbeh}{\mbox{\boldmath{$\hat{\beta}$}}}
\newcommand{\bSig}{\mbox{\boldmath{$\Sigma$}}}
\newcommand{\bLam}{\mbox{\boldmath{$\Lambda$}}}
\newcommand{\bGam}{\mbox{\boldmath{$\Gamma$}}}
\newcommand{\bOm}{\mbox{\boldmath{$\Omega$}}}
\newcommand{\bdel}{\mbox{\boldmath{$\delta$}}}
\newcommand{\cI}{\mbox{\boldmath{$\mathcal{I}$}}}
\newcommand{\cJ}{\mbox{\boldmath{$\mathcal{J}$}}}
\newcommand{\cY}{\mbox{\boldmath{$\mathcal{Y}$}}}
\newcommand{\xtil}{\tilde{\bX}}
\newcommand{\bbes}{\mbox{\boldmath{\scriptsize $\bbe$}}}
\newcommand{\bths}{\mbox{\boldmath{\scriptsize $\bth$}}}
\newcommand{\bX}{{\bf X}}
\newcommand{\bC}{{\bf C}}
\newcommand{\bB}{{\bf B}}
\newcommand{\bD}{{\bf D}}
\newcommand{\bG}{{\bf G}}
\newcommand{\bV}{{\bf V}}
\newcommand{\bH}{{\bf H}}
\newcommand{\bP}{{\bf P}}
\newcommand{\bQ}{{\bf Q}}
\newcommand{\bI}{{\bf I}}
\newcommand{\bS}{{\bf S}}
\newcommand{\bU}{{\bf U}}
\newcommand{\bfb}{{\bf b}}
\newcommand{\ba}{{\bf a}}
\newcommand{\bA}{{\bf A}}
\newcommand{\bZ}{{\bf Z}}
\newcommand{\bx}{{\bf x}}
\newcommand{\bm}{{\bf m}}
\newcommand{\bY}{{\bf Y}}
\newcommand{\bv}{{\bf v}}
\newcommand{\bR}{{\bf R}}
\newcommand{\bW}{{\bf W}}
\newcommand{\bw}{{\bf w}}
\newcommand{\bwun}{{\bf 1}}
\newcommand{\bzro}{{\bf 0}}
\newcommand{\pr}{\prime}
\newcommand{\half}{1/2}
\newcommand{\halffr}{\frac{1}{2}}
\newcommand{\Vh}{\bV^{\half}}
\newcommand{\Vmh}{\bV^{-\half}}
\newcommand{\frsm}[2]{\mbox{\small $\frac{#1}{#2}$}}
\newcommand{\nrm}[1]{\| #1 \|_\infty}
\newcommand{\Var}{\mbox{Var}}
\newcommand{\Cov}{\mbox{Cov}}
\newcommand{\cX}{\mathcal{X}}
\newcommand{\sumi}{\sum_{i=1}^n}
\newcommand{\Zbw}{\bar{Z}_w}
\newcommand{\intrl}{\int_{-\infty}^{\infty}}
\newcommand{\hB}{\hat{B}_n}
\newcommand{\tB}{\tilde{B}_n}
\newcommand{\chB}{\check{B}}
\newcommand{\cV}{\mathcal{V}}
\newcommand{\tUps}{\tilde{\Upsilon}}
\newcommand{\bi}{\begin{itemize}}
\newcommand{\ei}{\end{itemize}}
\newenvironment{Rinp}{
\begin{list}{}{
\setlength{\leftmargin}{0.5cm}}
 \item[]\begin{alltt}
}
{\end{alltt}\end{list}\normalsize\vspace{0.17cm}
}

\newcommand{\proof}[2]
{
\bigskip

\noindent {\bf {\small Proof. #1}} $\,$ {\footnotesize  #2} \hfill $\Box$

\bigskip
}
\newtheorem{thm}{Theorem}
\newtheorem{lemma}{Lemma}
\newtheorem{prop}{Proposition}

\begin{center}{\Large{\bf Subtleties in the interpretation of hazard ratios}}
 \end{center}
\vspace{0.2cm}

\begin{center}
{ \large Torben $\mbox{Martinussen}^1$, Stijn $\mbox{Vansteelandt}^2$ and Per Kragh $\mbox{Andersen}^1$}\\
 \vspace{2mm}
$\mbox{ }^1$ Department of Biostatistics\\
University of Copenhagen\\
\O ster Farimagsgade 5B, 1014 Copenhagen K, Denmark\\
 \vspace{2mm}
$\mbox{ }^2$ Department of Applied Mathematics, Computer Science and Statistics\\
Ghent University  \\
Krijgslaan 281, S9, B-9000 Gent, Belgium \\
\end{center}

\nothere{
\begin{center}
{ \large Torben Martinussen, Per Kragh Andersen}\\
Department of Biostatistics\\
University of Copenhagen\\
\O ster Farimagsgade 5B, 1014 Copenhagen K, Denmark\\
 \vspace{4mm}
{ \large Stijn Vansteelandt} \\
Department of Applied Mathematics, Computer Science and Statistics\\
Ghent University  \\
Krijgslaan 281, S9, B-9000 Gent, Belgium \\
\end{center}
}
%\newpage
\vspace{2cm}

\bigskip \setlength{\parindent}{0.3in} \setlength{\baselineskip}{24pt}

\centerline{\sc Summary}
\noindent

The hazard ratio is one of the most commonly reported measures of treatment effect in randomised trials, yet the source of much misinterpretation.   
This point was made clear by \cite{Hernan2010} in commentary, which
 emphasised that the hazard ratio contrasts populations of treated and untreated individuals who survived a given period of time, populations that 
 will typically fail to be comparable - even in a randomised trial - as a result of different pressures or intensities acting on both populations. The commentary has been very influential, but also a source of surprise and confusion. In this note, we aim to provide more insight into the subtle interpretation of hazard ratios and differences, by investigating in particular what can be learned about treatment effect from the hazard ratio becoming 1 after a certain period of time. Throughout, we will focus on the analysis of randomised experiments, but our results have immediate implications for the interpretation of hazard ratios in observational studies.
\noi
 \vspace{3mm}

\noi

\section{Introduction}

The popularity of the Cox regression model has contributed to the enormous success of the hazard ratio as a concise summary of the effect of a randomised treatment on a survival endpoint. Notwithstanding this, use of the hazard ratio has been criticised over recent years. \cite{Hernan2010} argued that selection effects render a causal interpretation of the hazard ratio difficult when treatment affects outcome. While the treated and untreated people are comparable by design at baseline, the treated people who survive a given time $t$ may then tend to be more frail (as a result of lower mortality if treatment is beneficial) than the untreated people who survive the given time $t$, so that the crucial comparability of both groups is lost at that time. \cite{AalenCook15} re-iterated Hern\'an's concern. They viewed the problem more as one of non-collapsibility \citep{MartVan_Collap}, which is a concern about the interpretation of the hazard ratio, though not about its justification as a causal contrast. In particular, they argued that the magnitude of the hazard ratio typically changes as one evaluates it in smaller subgroups of the population (e.g. frail people), even in the absence of interaction effects on the log hazard scale. In this paper, we aim to develop more insight into these matters. We will focus in particular on what can be learned about the treatment effect from the hazard ratio becoming 1 after a certain point in time. Throughout, we will assume that data are available from a randomised experiment on the effect of a dichotomous treatment $A$ (coded 1 for treatment and 0 for control) on a survival endpoint $T$, so that issues of confounding can be ignored, although our findings naturally extend to observational studies.   

%{\color{blue} 
\section{The Cox model and causal reasoning}
%}

Analyses of time-to-event endpoints in randomised experiments are commonly based on the Cox model
$$
\lambda(t;a)=\lambda_0(t)e^{\beta a},
$$
where $\lambda(t;a)$ denotes the hazard function of $T$ given $A=a$, evaluated at time $t$, and $\lambda_0(t)$ is the unspecified baseline hazard function.
%and let $\Lambda_0(t)=\int_0^t\lambda_0(s)\, ds$ denote the cumulative baseline hazard function.
This model implies that for all $t$
\begin{equation}\label{survival}P(T>t |A=1)=P(T>t |A=0)^{\exp (\beta)},\end{equation}
which suggests that the exponential of $\beta$ can be interpreted as
\[\exp(\beta) = \frac{\log P(T>t |A=1)}{\log P(T>t |A=0)}.\]
This represents a causal contrast (i.e., it compares the same population under different interventions). Indeed, let $T^a$ denote the potential event time we would see if the exposure $A$ is set to $a$. Then, since randomisation ensures that $T^a\cip A$ for $a=0,1$, we have that 
\[\exp(\beta) = \frac{\log P(T^1>t)}{\log P(T^0>t)}.\]
This shows that, under the proportional hazard assumption, $\exp(\beta)$ forms a relative contrast of what the log survival probability would be at an arbitrary time $t$ if everyone were treated, versus what it would be at that time if no one were treated. 

The log-transformation makes the above interpretation of $\exp(\beta)$, while causal, difficult. 
It is therefore more common to interpret $\exp(\beta)$ as a hazard ratio
\[\exp(\beta) = \frac{\lim_{h\rightarrow 0}P(t\le T< t+h|T\ge t,A=1)}{\lim_{h\rightarrow 0}P(t\le T< t+h|T\ge t,A=0)}= \frac{\lim_{h\rightarrow 0}P(t\le T^1< t+h|T^1\ge t)}{\lim_{h\rightarrow 0}P(t\le T^0< t+h|T^0\ge t)}.\]
Interpretation appears simpler now, but this is somewhat deceptive for two reasons. First, the righthand expression shows that $\exp(\beta)$ contrasts the hazard functions  with and without intervention for two separate groups of individuals, those who survive time $t>0$ with treatment ($T^1\ge t$) and those who survive time $t>0$ without treatment ($T^0\ge t$). Those groups will typically fail to be comparable if treatment affects outcome \citep{Hernan2010}. In particular, when treatment has a beneficial effect then, despite randomisation, the subgroup $T^1\geq t$ in the numerator will generally contain more frail people than the subgroup  $T^0\geq t$ in the denominator, where the frailest people may have died already. When viewed as a hazard ratio, $\exp(\beta)$ therefore does not represent a causal contrast. Second, the interpretation of $\exp(\beta)$ as a hazard ratio is further complicated by it being non-collapsible, so that its magnitude typically becomes more pronounced as one evaluates smaller subgroups of the study population \citep{MartVan_Collap,AalenCook15}. 

The summary so far is that, under the assumption of proportional hazards, the parameter $\exp(\beta)$ in the Cox model expresses a causal effect, namely the ratio of log survival probabilities with versus without treatment in the study population. Because this interpretation is not insightful, results are best communicated by visualising identity (\ref{survival}) in terms of estimated survival curves with versus without treatment \citep{Hernan2010}. This has the advantage that it provides better insight into the possible public health impact of the intervention, but the drawback that it does not permit a compact way of reporting and that survival curves do not provide an understanding of a possible dynamic treatment effect. To enable a more in-depth understanding, it is tempting to interpret  $\exp(\beta)$ as a hazard ratio, but then interpretation becomes subtle. 
We will demonstrate this in more detail in the next section, where we will investigate to what extent hazard ratios may provide insight into the dynamic nature of the treatment effect.
\nothere{
\footnote{I HAVE NOW DROPPED THE SECTION `FURTHER POINTS'. If the conditional HR is always below 1, then the conditional survival function is always higher in the treatment arm. This then also holds marginally, and thus we have a marginal HR below 1. So we can show this in general. That there is attenuation holds generally as well, I think. We should first see if there is such general result before restricting to a specific model.
{\color{blue}TM: I have spend some time on showing this in general, but so far without success. I've also searched the literature, and did not find a general result  for the Cox model. There is a result for logistic regression: Neuhaus, Kalbfl and Hauck (1991), but they only show it with a random int (corresp to frality in our setting), and random coef.}
}
}

\section{Time-varying hazard ratios are not causally interpretable}

Consider now a study where the hazard ratio changes with time in the following sense: 
\begin{equation}
\label{Cox_CP}
\frac{\lambda(t;A=1)}{\lambda(t;A=0)}= \left\{ \begin{array}{ll}
\exp(\beta_1) & \textrm{if $t\le \nu$}\\
\exp(\beta_2) & \textrm{if $t> \nu$}\\
\end{array} \right.
\end{equation}
with $\beta_1\neq \beta_2$, where $\nu>0$ denotes the change point, which we assume to be known based on subject matter knowledge. Suppose in particular that $\beta_1<0$ and $\beta_2=0$. 
This is commonly interpreted as if treatment is initially beneficial, but becomes ineffective from time $\nu$ onwards. 
 We present a practical example in Section 4.2 where this situation arises.
In the following sections, we will reason whether such interpretation is justified, and thus whether hazard ratios permit a dynamic understanding of the treatment effect.%\footnote{IT WOULD BE NICE IF WE HAD DATA FROM A REAL STUDY WHERE THE HAZARD RATIO INDEED BECAME 1 AFTER SOME TIME AND/OR EXAMPLES FROM REAL STUDIES WHERE SUCH INTERPRETATION WAS MADE.
%{\color{blue}TM: I have included a second application considering this aspect.}
%}

\subsection{A closer look at the causal mechanism}

To develop a greater understanding, we will first develop insight into data-generating processes (DGP) that could give rise to (\ref{Cox_CP}).
Let $Z$ represent the participants' unmeasured baseline frailty (higher means more frail), which affects $T$, but is independent of $A$ by randomisation.
Suppose that the hazard function $\lambda(t;a,z)$ of $T$ given $A=a$ and $Z=z$ satisfies 
\begin{equation}\label{Coxz}
\lambda(t;a,z)=z\lambda^*(t;a),\end{equation} for some function $\lambda^*(t;a)$, and let 
$Z$ be Gamma distributed with mean 1 and variance $\theta$.
%\sim \Gamma (1/\theta,\theta)$ with $\theta=1$, i.e. $Z$ is exponentially distributed with mean 1. T
This specific choice is not essential, however, and we later consider a situation where $Z$ is  binary. 
%We will moreover show that this model does not contradict model (\ref{Cox_CP}). In this section, we will then evaluate what the above model (\ref{Cox_CP}) implies among individuals that share the same $Z$.
  
We will investigate what choices of $\lambda^*(t;a)$ give rise to model (\ref{Cox_CP}).
With $\phi_Z(u)=E(e^{-Zu})$ the Laplace transform associated with the distribution of $Z$, the following relationship between the hazard function of interest $\lambda(t;A)$, and $\lambda^*(t;A)$, can be shown to hold:
$$
\Lambda^*(t;a)=\phi_Z^{-1}(e^{-\Lambda(t;a)})=\frac{1-e^{-\theta\Lambda(t;a)}}{\theta e^{-\theta\Lambda(t;a)}},
$$
where $\Lambda(t;a)=\int _0^t \lambda(s;a)\, ds$ and similarly with $\Lambda^*(t;a)$.
%Since we have decided on the structure of $\Lambda(t;a)$, the Cox model with a change point, we can find $\Lambda^*(t;a)$ by a direct calculation
%using the expression in the latter display, and then also express $\lambda(t;a,z)$ via  $\lambda(t;a,z)=z\lambda^*(t;a)$. 
Simple calculations then show that model (\ref{Coxz}) implies model (\ref{Cox_CP}) when
\begin{equation}
\label{Cox_AZ}
\lambda(t;A,Z)= \left\{ \begin{array}{ll}
Z\lambda_0(t)e^{\beta_1A}\exp{\{\theta\Lambda_0(t)e^{\beta_1A}\}} & \textrm{if $t\le \nu$}\\
 Z\lambda_0(t)e^{\beta_2A}\exp{\left \{\theta\Lambda_0(\nu)e^{\beta_1A}+\theta\Lambda_0(\nu,t)e^{\beta_2A}\right \}}& \textrm{if $t> \nu$}\\
\end{array} \right.
\end{equation}
where $\Lambda_0(\nu,t)=\int_{\nu}^t\lambda_0(s)\, ds$.  For subjects with given $Z=z$, it follows that the conditional hazard ratio, which we term $\mbox{HR}_{Z}(t)$, is
$$\frac{\lambda(t;A=1,Z)}{\lambda(t;A=0,Z)}= \left\{ \begin{array}{ll}
\exp(\beta_1)\exp{[\theta\Lambda_0(t)\{\exp(\beta_1)-1\}]} & \textrm{if $t\le \nu$}\\
\exp(\beta_2)\exp{[\theta\Lambda_0(\nu)\{\exp(\beta_1)-1\}+\theta\Lambda_0(\nu,t)\{\exp(\beta_2)-1\}]}& \textrm{if $t> \nu$}\\
\end{array} \right.
$$
For $\beta_2=0$, this simplifies to
$$\frac{\lambda(t;A=1,Z)}{\lambda(t;A=0,Z)}= \left\{ \begin{array}{ll}
\exp(\beta_1)\exp{[\theta\Lambda_0(t)\{\exp(\beta_1)-1\}]} & \textrm{if $t\le \nu$}\\
\exp{[\theta\Lambda_0(\nu)\{\exp(\beta_1)-1\}]}& \textrm{if $t> \nu$}\\
\end{array} \right.
$$
which is smaller than 1 at all time points when, as previously assumed, $\beta_1<0$. We conclude that treatment appears beneficial at all times
amongst individuals with the same value $z$ of $Z$, regardless of $z$. This contradicts the earlier, na\"{\i}ve interpretation that, across all individuals combined, treatment is ineffective from time $\nu$ onwards. 

The root cause of these contradictory conclusions is the fact that the hazard ratio at a given time does not express a causal effect
% as also explained by 
 \citep{Hernan2010}.
 This has nothing to do with model misspecification, as all considered models hold by construction.
% but not model misspecification (as all considered models were considered to hold).
In particular, when $\beta_2=0$ and $\theta=1$,  then %({\color{blue} Torben: Check below expression})
$$
\frac{E(Z|T>t,A=1)}{E(Z|T>t,A=0)}= \left\{ \begin{array}{ll}
\exp{\{\Lambda_0(t)(1-e^{\beta_1})\}} & \textrm{if $t\le \nu$}\\

\frac{\exp{\{\Lambda_0(\nu)\}} +\Lambda_0(\nu,t)}{\exp{\{\Lambda_0(\nu)e^{\beta_1}\}} +\Lambda_0(\nu,t)}  & \textrm{if $t>\nu$}.
\end{array} \right.
$$
While $Z$ is independent of $A$ by randomisation, it is thus no longer so amongst subgroups of survivors, where  we are left with more frail subjects in the active treatment group: $E(Z|T>t,A=1)>E(Z|T>t,A=0)$. That selection takes place does not rely on $Z$ being Gamma-distributed; see  Appendix A.1, where we consider a situation with a binary frailty variable. 

To illustrate further, we generated a single data set with
$Z$ binary with $P(Z=0.2)=0.2$ and $P(Z=1.2)=0.8$ corresponding to low and high risk groups. The treatment variable $A$ was  binary with $P(A=1|Z)=0.5$. Further,  
we took $\beta_1=-\log{(2)}$,
$\beta_2=0$, the baseline hazard function  $\lambda_0(t)=0.4$, and the change point  $\nu=4$. We took a large sample size $n=20000$ so that sampling variability is small, and induced censoring according to a uniform distribution on $[0,10]$, 
randomly for half of the individuals, and the rest censored at $t=8$, resulting in an overall censoring percent of approximately 19\%. 
The result from a change point Cox analysis was $\hat\beta_1=-0.67$ (SE 0.018) and $\hat\beta_2=-0.03$ (SE 0.034). 
The corresponding hazard ratios are $0.51$ for $t\le 4$ and $1.03$ for $t>4$. In contrast, the conditional hazard ratio in Figure 1
is seen to be smaller than 1 at all times, not only in the first initial period. 
Indeed, no matter the value of $Z$, it lies between 0.33 and 0.80 in the considered time period from 0 to 8.

That we are misled by the hazard ratios calculated from the extended
 Cox regression analysis is due to selection taking place. This is shown in Figure 2, where 
 we have  plotted $E(Z|T>t,A=a)$ for $a=0,1$.  Since the treatment has a beneficial effect, we are left with more and more frail subjects  in this group compared to the untreated group. In Figure 2, the blue curve ($A=1$) therefore lies consistently above the green curve ($A=0$). 
  
 %  {\color{blue}
As just shown,
the DGP  given by $\lambda(t;A,Z)$ with $\nu=\infty$ implies  the marginal Cox model.
It is tempting to interpret $\mbox{HR}_{Z}(t)$ as a causal hazard ratio, but this only holds under further untestable assumptions as shown in Section 3.2 and 3.4.
In  Appendix A.2, we formulate a more general DGP so that both the marginal Cox model and the model $\lambda(t;A,Z)$ are correctly specified.
 %}
 %Since we assume
% the true DGP is governed by $\lambda(t;A,Z)$ this essentailly means that $Z$ captures every aspects about the risk of dying except for a potential effect of the treatment. 
%We write 
%$$\lambda(t;A)=\{\exp(\beta_1A)I(t\leq \nu)+ \exp(\beta_2A)I(t>\nu)\} \lambda_0(t), 
%$$
%The hazard function  $ \lambda_0(t)$ is not further specified. 
%Note that  the Cox model is obtained by taking $\nu$ equal  to infinity, leading to only one HR, $e^{\beta_1}$.

%We use the potential outcome notation (explain here) $P(T^a>t)$ and 
%write the corresponding hazard function as $\lambda(t;A)$, and 

\subsection{Towards a causal hazard ratio}

To remedy this problem, it seems intuitively of interest to evaluate conditional hazard ratios
\[\frac{\lim_{h\rightarrow 0}P(t\le T< t+h|T\ge t,A=1,Z=z)}{\lim_{h\rightarrow 0}P(t\le T< t+h|T\ge t,A=0,Z=z)}= \frac{\lim_{h\rightarrow 0}P(t\le T^1< t+h|T^1\ge t,Z=z)}{\lim_{h\rightarrow 0}P(t\le T^0< t+h|T^0\ge t,Z=z)},\]
for a large collection of baseline variables $Z$, such that those who survive time $t>0$ with treatment ($T^1\ge t$) are comparable to those who survive time $t>0$ without treatment ($T^0\ge t$), but have the same covariate values $z$. 
%{\color{blue}
 Such comparability would be attained if %when $Z$ deterministically predicts $T^1$ and $T^0$ such that
 \begin{equation}
 \label{cip}
 T^1\cip T^0|Z. 
\end{equation}
%}
%Less demandingly, we will here assume that $\{I(t\le T^1< t+h),I(T^1\ge t)\}\cip \{I(t\le T^0< t+h),I(T^0\ge t)\}|Z$. 
\noindent
Under this assumption, it follows via Bayes' rule that the righthand side of the above identity equals
\begin{equation}\label{causalHR}\frac{\lim_{h\rightarrow 0}P(t\le T^1< t+h|T^0\ge t,T^1\ge t,Z=z)}{\lim_{h\rightarrow 0}P(t\le T^0< t+h|T^0\ge t,T^1\ge t,Z=z)}.\end{equation}
This estimand expresses the instantaneous risk at time $t$ on treatment versus control for the principal stratum of individuals with covariates $z$ who would have survived up to time $t$, no matter what treatment. It represents a causal contrast, which is closely related to the so-called survivor average causal effect \citep{10.2307/2669382}. We will refer to it as the conditional causal hazard ratio.

Unfortunately, the assumption \eqref{cip} is untestable and biologically implausible, as it is essentially impossible to believe that one can get hold of all predictors of the event time such that knowledge of the event time without treatment does not further predict the event time with treatment. Furthermore, even if one could get hold of all such predictors, then because $Z$ would probably
carry so much information about the event time, one would logically expect the numerator and denominator of (\ref{causalHR}) to be so close to 0 or 1 that it would render the conditional causal hazard ratio essentially meaningless.
%\footnote{Your additional calculations intended to investigate if we include more and more $Z$'s are we then approaching a situation where $T^0$ and $T^1$ are cond. indep., could perhaps be placed in the appendix and linked to this paragraph.}
Below, we will therefore focus on the marginal causal hazard ratio $\mathrm{HR}(t)$
obtained from \eqref{causalHR} with $Z$ empty:
%{\color{blue}
 \begin{equation}
 \label{HRt}
\mathrm{HR}(t)=\frac{\lim_{h\rightarrow 0}P(t\le T^1< t+h|T^0\ge t, T^1\ge t)}{\lim_{h\rightarrow 0}P(t\le T^0< t+h|T^0\ge t,T^1\ge t)}.
 \end{equation}
%}
%, which is defined like $\mathrm{HR}(t|z)$, but with $Z$ empty.

\subsection{Why causal hazard ratios are not identified without strong assumptions}

The reason why the causal hazard ratio \eqref{HRt} is not identifiable without invoking strong assumptions is that it attempts to answer an overly  ambitious question. Imagine a trial that randomises participants over an implanted medical device (e.g. a stent or a pacemaker) versus no treatment (or placebo). Suppose that the medical device gradually deteriorates and
stops being operational after some time $\nu$. Then we would say that treatment no longer works from time $\nu$ onwards. This would correspond with $\mathrm{HR}(t)=1$ for $t\ge \nu$.
However, how could data from a randomised trial be informative about the effect of treatment after time $\nu$ when no information is collected on the times at which the medical device is operational or not? To learn about the treatment effect at each time $t$, we should ideally need data $A_t$ on whether ($A_t=1$) or not ($A_t=0)$ the device is operational at that time. When the operation time is ignorable, then one may learn about the treatment effect at each time $t$ through contrasts of the form
\[\frac{\lim_{h\rightarrow 0}P(s\le T< s+h|T\ge s,\overline{A}_{t-}=1,A_t=1)}{\lim_{h\rightarrow 0}P(s\le T< s+h|T\ge s,\overline{A}_{t-}=1,A_t=0)},\]
for all $s\geq t$, where $\overline{A}_{t-}$ is the information generated by all $A_u$, $u<t$.
It is unsurprising that without detailed data on the operation times of each device, strong assumptions are needed to develop insight into the dynamic nature of the treatment effect. Likewise, consider a trial that randomises participants over a once-daily treatment regimen versus placebo. Then we would say that treatment no longer works from time $\nu$ onwards when, from that time onwards, patients with the same history of treatment experience the same outcomes (in distribution), whether or not they continue treatment. To infer when treatment becomes ineffective, a multi-stage design is ideally needed where patients on the treatment arm may randomly be switched to the control arm at designated points in time. Without such design, some progress can still be made with data on daily pill intake. However, without such design and data, inferring when treatment becomes ineffective remains an ambitious undertaking.

\subsection{Sensitivity analysis}

Some further insight into the magnitude of the causal hazard ratio at a given time $t$ can be obtained under the monotonicity assumption that no one is harmed by treatment, i.e.
\[T^1\geq T^0 \quad \mathrm{ with \ probability \ 1}.\]
Under this assumption, we can write
\[P(t\le T< t+h|T\ge t,A=0)=P(t\le T^0< t+h|T^0\ge t)=P(t\le T^0< t+h|T^0\ge t,T^1\ge t)\]
and
\begin{eqnarray*}
P(t\le T< t+h|T\ge t,A=1)&=&P(t\le T^1< t+h|T^1\ge t)\\
&=&P(t\le T^1< t+h|T^0\ge t,T^1\ge t)\pi(t)\\
&&+P(t\le T^1< t+h|T^0< t,T^1\ge t)\left\{1-\pi(t)\right\},
\end{eqnarray*}
with 
\begin{eqnarray*}
\pi(t)&\equiv& P(T^0\ge t|T^1\ge t)=\frac{P(T^0\ge t)}{P(T^1\ge t)}=\frac{P(T\ge t|A=0)}{P(T\ge t|A=1)}.
\end{eqnarray*}
It follows that 
\begin{equation}\label{selection}
\frac{\lambda(t;A=1)}{\lambda(t;A=0)}=\mathrm{HR}(t)\left[\pi(t)+\mathrm{SR}(t)\left\{1-\pi(t)\right\}\right],\end{equation}
where
\[\mathrm{SR}(t)\equiv \frac{\lim_{h\rightarrow 0}P(t\le T^1< t+h|T^0< t,T^1\ge t)}{\lim_{h\rightarrow 0}P(t\le T^1< t+h|T^0\ge t,T^1\ge t)}.\]
compares the instantaneous risk at time $t$ on treatment for individuals who remained event-free at time $t$ thanks to treatment versus individuals who would have been event-free at time $t$ regardless of treatment. It can be viewed as a selection effect and is therefore termed a sensitivity ratio (SR). Expression (\ref{selection}) thus conveys that the dependence of the hazard ratio ${\lambda(t;A=1)}/{\lambda(t;A=0)}$ on time $t$ may differ from the causal hazard ratio on time as a result of selection effects. In particular, when the observed hazard ratio ${\lambda(t;A=1)}/{\lambda(t;A=0)}$ is constant, e.g. suggesting a beneficial treatment effect of 0.8, then this will often correspond with a more pronounced treatment effect $\mbox{HR}(t)\le 0.8$. 
This results from \eqref{selection}, if $\mbox{SR}(t)\ge 1$, as $0 \le \pi(t)\le 1$.
The causal hazard ratio $\mbox{HR}(t)$ then need not be constant.
%changing with time 
This may be the result of varying treatment effectiveness over time, but may also be the result of the patient population (the principal stratum $T^0\ge t,T^1\ge t$)  changing with time.
It follows from the above results that 
\[\mathrm{HR}(t)=\frac{\lambda(t;A=1)}{\lambda(t;A=0)}\left[\pi(t)+\mathrm{SR}(t)\left\{1-\pi(t)\right\}\right]^{-1}.\]
All terms in the righthand side, apart from $\mbox{SR}(t)$, are identified from the observed data.
The above expression may therefore be used as the basis of a sensitivity analysis where the user tries different choices of $\mbox{SR}(t)$. However, it is not clear what would be reasonable values of $\mbox{SR}(t)$, and the monotonicity assumption may also not be plausible. 

\nothere{
Let $S_t(a)=I(T^a>t)$, $a=0,1$,
and let  $Z$ denote some (potentially un-measured) variable so that, for $u>t$, 
\begin{equation}
\label{Ass_add}
P\{T^a>u|S_t(a)=1, S_t(1-a),Z=z\}=P\{T^a>u|S_t(a)=1,Z=z\},
\end{equation}
Hence, 
if one had observed $Z=z$ for an individual alive at time $t$ under exposure $a$ , 
further knowing whether the individual would also be alive at time $t$ under exposure $1-a$ does not further improve our ability to predict the survival status $T=T^{a}$   at later timepoints observed under exposure $a$.
%given $V$, $U$, and $S_t(a)=1$, the staus of the other counterfactual 
Note that \eqref{Ass_add} holds if $T^0$ and $T^1$ are conditionally independent given $Z$.}

%{\color{blue}
Assume instead that 
there is a $Z$ so that \eqref{cip} holds, and 
so that the DGP is governed by
\begin{equation}
\label{GenF}
\lambda(t;A,Z)=\exp{\{\psi_0(t,A)+\psi_1(t,Z)\}}
\end{equation}
for some general functions $\psi_0$ and   $\psi_1$, 
then
$$
\mbox{HR}(t)=\frac{E\{P(T^{1}=t|T^{1}\ge t, Z)| T^0\ge t, T^1\ge t\}}{E\{P(T^{0}=t|T^{0}\ge t, Z)| T^0\ge t, T^1\ge t\}}=
\exp{\{{\psi_0(t,1)-\psi_0(t,0)\}}=\mbox{HR}_{Z}(t)
}
$$
Thus, if  \eqref{cip} and \eqref{GenF} hold  then  $\mbox{HR}(t)=\mbox{HR}_Z(t)$.
%and if $Z$ is the empty set (marginal independence) then the Cox HR is equal to $\mbox{HR}(t)$.
%We illustrate this in Section \ref{MRC RE01 study}.
%e.g. different choices of $\gamma\geq 0$ in the user-specified model $\mathrm{SR}(t)=\exp(-\gamma t)$. Here, $\gamma=0$ when $T^1\cip T^0$, which is generally implausible. 
In the worked application in Section 4.1 we compare $\mbox{HR}(t)$ to the Cox hazard ratio while modelling the correlation between $T^0$ and $T^1$ letting $Z$ be Gamma distributed with varying variance so that \eqref{cip} and \eqref{GenF} hold.
%{\color{blue}
As a further illustration we now describe how to simulate data so that 
\eqref{cip} and \eqref{GenF} are fulfilled, and so that the marginal Cox model, conditioning only on $A$, is also correctly specified.
Take $A$, $Z$, $V_0$ and $V_1$ to be independent so that 
$Z$ is Gamma distributed with mean 1 and variance $\theta$, $V_0$ and $V_1$ are exponentially distributed with mean 1, and the  exposure $A$ is binary with $P(A=1)=1/2$. Then let
$T^0=\frac{1}{\theta}\log{(\frac{\theta}{Z}V_0+1)},$
$T^1=\frac{1}{\theta e^{\beta}}\log{(\frac{\theta}{Z}V_1+1)},$ and let
$T=(1-A)T^0+AT^1.$
It follows directly that \eqref{cip}  holds, and further 
 that the marginal Cox model is correctly specified,
  and also 
that
$$
\lambda(t;A,Z)=Z\exp{\{\beta A+\theta e^{\beta A}t\}}
$$
so \eqref{GenF} also holds.
Thus, $\mbox{HR}(t)=\mbox{HR}_{Z}(t)$,
and
$$
\mbox{HR}(t)=e^{\beta}\exp{\{t\theta(e^{\beta}-1)\}}.
$$
% compare with expression on top of page 6 in paper (I have taken $\lambda_0(t)=1$ here,
%and therefore $\Lambda_0(t)$ is equal to $t$). This shows that we can generate data so that (a) and (b1) are true.
Different values of $\theta$ will give different values of Kendall's $\tau$  corresponding to different correlation between $T^0$ and $T^1$. 
%{\color{blue}
In the specific setting,
we have $\tau=\theta/(\theta+2 )$.
%}
Figure 3 displays $\mbox{HR}(t)$ in scenarios with different values of  Kendall's $\tau$.
 Note that the marginal Cox model induces the hazard ratio $e^{\beta}$, which is taken to 
be $e^{\beta}=0.5$. It is seen that 
$\mbox{HR}(t)$ is equal to $e^{\beta}$ only under independence between $T^0$ and $T^1$; otherwise $\mbox{HR}(t)$ looks like a decreasing function in $t$ with 
 $\mbox{HR}(0)=e^{\beta}$. When  $T^0$ and $T^1$ are highly correlated (Kendall's $\tau=0.93$) then  $\mbox{HR}(t)$ is sharply decreasing towards zero. In Appendix A.4 we describe another DGP where again \eqref{cip} holds and the marginal Cox model is correctly specified, but  $\mbox{HR}(t)$ and $\mbox{HR}_{Z}(t)$ are different due to failure of \eqref{GenF}.

\nothere{In a hope to lessen the degree of sensitivity, one may also consider estimation of the conditional causal hazard ratio
\[\mathrm{HR}(t|z)\equiv \frac{\lim_{h\rightarrow 0}P(t\le T^1< t+h|T^0\ge t,T^1\ge t,Z=z)}{\lim_{h\rightarrow 0}P(t\le T^0< t+h|T^0\ge t,T^1\ge t,Z=z)},\]
by identifying it as
\[\mathrm{HR}(t|z)=\frac{\lambda(t;A=1,Z=z)}{\lambda(t;A=0,Z=z)}\left[\pi(t)+\left\{\mathrm{SR}(t)-1\right\}\left\{1-\pi(t)\right\}\\right]^{-1},\]}
\nothere{where 
\[\mathrm{SR}(t|z)\equiv \frac{\lim_{h\rightarrow 0}P(t\le T^1< t+h|T^0< t,T^1\ge t,Z=z)}{\lim_{h\rightarrow 0}P(t\le T^1< t+h|T^0\ge t,T^1\ge t,Z=z)}.\]
Conditioning on strong prognostic factors $Z$ of the event time may render the assumption that $\mathrm{SE}(t|z)=1$ more plausible. However, caution is warranted because non-collapsibility of the hazard ratio also increases the plausibility of $\mathrm{SR}(t|z)$-values far from 1, and moreover, we generally expect the assumption that $\mathrm{SR}(t|z)=1$ to be implausible.}

\subsection{Hazard differences}

Under some conditions, hazard differences - as opposed to hazard ratios - have a more appealing interpretation, apart from them being collapsible \citep{MartVan_Collap}. In particular, suppose that the additive hazards model
\begin{equation}\label{ah}
\lambda(t;A,Z)=\psi(t)A+\omega(t,Z),
\end{equation}
holds for general functions $\psi$, $\omega$ and baseline covariates $Z$.
%{\color{blue} 
It has been shown that under this model, the baseline exchangeability of treated and untreated individuals w.r.t. the covariates $Z$, as guaranteed by randomisation, 
implies  that  $A\cip Z| T>t$ \citep{Vansteelandt:2014aa,AalenCook15}. 
%When this model holds, selection effects affecting the distribution of $Z$ are thus ruled out. 
However,
since this is equivalent to 
\begin{equation}
\label{ah_cip}
P(Z=z\cond T^0>t)=P(Z=z\cond T^1>t),
\end{equation}
it does not mean that there is a balance in the risk set, because \eqref{ah_cip} shows that this is  a comparison between two different groups of people: those with $T^0>t$ and those with $T^1>t$.
%Unfortunately, this does not suffice for 
Indeed, the hazard difference
\[\psi(t)=\lim_{h\rightarrow 0}P(t\le T^1< t+h|T^1\ge t)-\lim_{h\rightarrow 0}P(t\le T^0< t+h|T^0\ge t)\]
may not represent a causal contrast. 
%}
%{\color{blue}
If \eqref{cip} holds, e.g. that 
$Z$ is  a sufficiently rich collection of variables that includes $T^0$  and $T^1$ (or deterministically predicts $T^0$ and $T^1$)
then
%For this, one needs $Z$ to be a sufficiently rich collection of variables that includes $T^0$ (or deterministically predicts $T^0$). In that case, 
the above hazard difference reduces to 
\[ \psi(t)=\lim_{h\rightarrow 0}P(t\le T^1< t+h|T^0\ge t,T^1\ge t)-\lim_{h\rightarrow 0}P(t\le T^0< t+h|T^0\ge t, T^1\ge t),\]
which represents a causal contrast. This is shown in Appendix A.5.
%}
While it is implausible
to have such rich collection of data that it essentially deterministically predicts $T^0$, interestingly,
 when model (\ref{ah}) holds for $Z$ including $T^0$, then it can be fitted without data on $T^0$. Indeed, by 
collapsibility of the hazard difference \citep{MartVan_Collap}, the hazard difference $\psi(t)$ can be consistently estimated via Aalen least squares estimation in the unadjusted model
\[\lambda(t;A)=\psi(t)A.\]
Unfortunately, however, the additive structure of model (\ref{ah}) will often be unlikely satisfied w.r.t. a rich collection of variables that predict $T^0$, and thus the practical implications of the above reasoning remain limited. For instance, in the simulation study of Section 3.1, the marginal Aalen additive hazards model fits the data perfectly, because $A$ is binary, but also suggests a beneficial effect of the treatment in the first 4 years, which then disappears, see Figure 4. 

 %The reason being that the exchangeability due to randomization only holds at time $t=0$, it is lost when $t>0$. 
 \nothere{If the study runs in an interval $[0,\tau]$ with $\tau$ not so much larger than $\nu$ one may even choose parameters so that $\mbox{HR}_Z(t)<1$ for all $t\in [0,\tau]$ but with $\beta_2>0$ so one might wrongly conclude that the treatment is harmful after time point $\nu$, although arguing based on $\mbox{HR}_Z(t)$, it seems to be  beneficial for all $t\in [0,\tau]$. }
 %The reason why we are cheated is as mentioned due to that exchangeability is lost for $t>0$ due to 
 
\nothere{One may alternatively work on the hazard difference scale and write the hazard difference 
\[\lambda(t;A=1,Z=z)-\lambda(t;A=0,Z=z)\]
as \[\mathrm{HD}(t|z)+\mathrm{SD}(t|z)\frac{P(T\ge t|A=0,Z=z)-P(T\ge t|A=1,Z=z)}{P(T\ge t|A=0,Z=z)},\]
where
\[\mathrm{HD}(t|z)\equiv \lim_{h\rightarrow 0}P(t\le T^1< t+h|T^0\ge t,T^1\ge t,Z=z)-P(t\le T^0< t+h|T^0\ge t,T^1\ge t,Z=z),\]
and
\[\mathrm{SD}(t|z)\equiv \lim_{h\rightarrow 0}P(t\le T^1< t+h|T^0< t,T^1\ge t,Z=z)-P(t\le T^1< t+h|T^0\ge t,T^1\ge t,Z=z).\]
Let $Z$ include a large number of possibly unmeasured prognostic factors of $T$ (possibly depending on time) so that 
individuals who would or would not survive time $t$ without treatment, but have the same covariate data $Z$ at time $t$ and would all survive time $t$ on treatment, are comparable (in terms of survival on treatment). Then $\mathrm{SD}(t|z)=0$. 
The conditional hazard difference then equals the causal hazard difference, and can be identified by collapsibility of the hazard difference when $A$ and $Z$ exert additive effects on the hazard.}

\section{Applications}
\subsection{MRC RE01 study}
\label{MRC RE01 study}

As an illustration, we reconsider the kidney cancer data described in \cite{WhiteRoyston_09}. These data
are from the MRC RE01 study that was a randomised
controlled trial comparing  interferon-$\alpha$ (IFN) treatment 
with the best supportive care and hormone treatment with
medroxyprogesterone acetate (control) in patients with metastatic
renal carcinoma. We use the same 347 patients as in White and Royston
(2009). In this illustrative analysis we  consider only the first 30 months of follow up. The median follow-up time was 242 days, and 85\% of the patients died within the considered time frame.
The two Kaplan-Meier estimates in Figure \ref{Fig: KM} contain all the available information about the treatment effect.
IFN treatment seems superior to the standard treatment (control), although a supremum test comparing the two survival curves results in a non-significant p-value of 0.09.
The score process plot of \cite{lin:wei:ying:1993} was calculated using the R-package {\tt timereg}  \citep{ms06}, giving no evidence against the proportional hazards assumption (P 0.43 according to the supremum test). We therefore fitted the Cox model, giving the estimate 
$\hat\beta=-0.29$ (SE 0.12, P 0.01), showing evidence that the IFN treatment reduces the risk of dying for these patients. The hazard ratio of 
0.75 expresses the magnitude of the causal effect, but interpretation is subtle as explained before.
As a sensitivity analysis we now assume  \eqref{cip} and \eqref{GenF}, and take $Z$ to be Gamma distributed with mean 1 and variance $\theta$. 
This model fits the observed data equally well and thus cannot be refuted based on the observed data. 
%{\color{blue}
The parameter $\theta$ expresses the correlation between $T^0$ and $T^1$ and was chosen to give a Kendall's $\tau$ of 0.1 and 0.2 corresponding to estimated Kendall's
$\tau$ concerning  the correlation between lifetimes for dizygotic and monozygotic Danish twins \citep{Scheike:2015aa}.
A  higher correlation between $T^0$ and $T^1$, corresponding to  a Kendall's $\tau$ of 0.3, was also considered.
%} 
Figure 6 displays $\mbox{HR}(t)$ under these three  scenarios. It is seen that $\mbox{HR}(t)$ is smaller than 0.75 at all times and decreases with time indicating a stronger treatment effect
 when comparing the instantaneous risk at time $t$ on treatment versus control for the principal stratum of individuals  who would have survived up to time $t$, no matter what treatment. This is more pronounced with the larger correlation.
 In view of these subtleties of interpretation, in this section, we will focus on a number alternative ways of describing the treatment effect.

\nothere{
We will therefore first consider alternative effect measures that have a causal interpretation.
From the probabilistic index model interpretation\footnote{NEEDS MORE EXPLANATION, BOTH WHAT THIS MODEL IS AND HOW YOU GET TO 1.34.} of the Cox model we may conclude that, on average, the odds of living longer, when comparing those randomized to  IFN with those randomized to the control treatment, is estimated to 
1.34 (95\% CI, 1.06 to 1.68). 
}
We may alternatively use the Cox model to estimate the relative risk function 
\[\mbox{RR}(t)=\frac{P(T\le t|A=1)}{P(T\le t|A=0)}=\frac{P(T^1\le t)}{P(T^0\le t)},\] 
which can be estimated consistently by
$$
\widehat{\mbox{RR}}(t)=
\frac{1-\exp{\{-\hat \Lambda_0(t)e^{\hat\beta}\}}}{1-\exp{\{-\hat \Lambda_0(t)\}}},
$$
where $\hat\beta$ is the Cox partial likelihood estimator and $\hat\Lambda_0(t)$ is the corresponding Breslow estimator. 
% In the case where there are additional covariates, $V$, the Cox model may still be used to estimate the relative risk function but this estimate will depend on the specific value of $V$, which may be inconvenient. We elaborate on this point in the Concluding remarks section. Using large sample results for the Cox model one can calculate 95\% pointwise as well as uniform confidence bands for 
%$\mbox{RR}(t)$. This is illustrated in Section 8, where we study data on kidney cancer, specifically see Figure 6.
This estimate, along with 95\% pointwise and uniform confidence bands, is displayed in Figure 7, where we see that the estimated relative risk function is below 1 at all times, in favour of the IFN treatment. For instance, it is seen that the relative risk at one year is estimated to be approximately 0.85, and, judging from the 95\% confidence bands (dashed curves), this is close to being significant. A uniform test over the considered time span is also close to being significant judging from the  95\% uniform confidence bands.

Another way of quantifying the treatment effect is by using the restricted mean survival time (RMST) \citep{Uno14,Zhao_RMST}.
The RMST up to time $t$ is defined as $\mbox{RMST}(t) = E\{\mbox{min}(T, t)\}$.
This is
the area under the survival curve of $T$ up to time $t$ and can easily be estimated using the corresponding Kaplan-Meier  curve up to time $t$ (see  \cite{Zhao_RMST} for details on inference). Contrasts of the $\mbox{RMST}(t)$ corresponding to different (randomised) treatment groups therefore carry a causal interpretation.
The restricted mean time lost, $\mbox{RMTL}(t)$ is defined as $t-\mbox{RMST}(t)$.
%Specifically, we calculated  $\mbox{RMTL}(30)$ under the two treatments.
Here, the `months of life lost up to 30 months' is given by $\mbox{RMTL}(30)$ and estimated to 17.3 for the IFN treatment and to 19.8
 for the control treatment. The ratio of these two (IFN vs control) is 0.87 (95\% CI, 0.70 to 1.04).  Thus, on IFN treatment there is a 13\% less loss of lifetime compared to the control treatment during the first  30 months of follow up.

We also fitted Aalen's additive hazard model
\begin{equation}
\label{AalenM}
\lambda(t;A,V)=\beta_0(t)+\psi(t)A+\beta_1(t)^TV,
\end{equation}
where $V$ includes days from metastasis to randomization (log-transformed), 
WHO performance status (0, 1 and 2; with group 0 and 1 collapsed into one group), and 
Haemoglobin (g/dl).
%{\color{blue} 
 As the treatment variable $A$ is independent of the other covariates, the above Aalen model is collapsible meaning that the interpretation of  $\psi(t)$ is the same in the conditional and marginal model. The Cox model does not have this property.
%}
Model \eqref{AalenM} appeared to fit the data well, using the tools described in Chapter 5 in  \cite{ms06}; specifically no interaction between the treatment indicator and the baseline risk factors was found.
We estimated $\hat{\Psi}(t)=\int_0^t\psi(s)ds$  both from the conditional model and the marginal model, and the two estimators were  almost identical,
which, as pointed out,  should be the case if model \eqref{AalenM} is correctly specified.
%Hence we applied the following model
Using model \eqref{AalenM}, we next tested the null hypothesis $\psi(t)=\psi$ of a constant effect (P 0.58), which was subsequently estimated to be $\hat\psi=-0.02$ (SE, 0.009).  If $T^0\cip T^1|V$ (or if the addition of additional variables $Z$ conditional on which $T^0$ and $T^1$ become independent, does not change the additive structure of the model), then this is also the causal hazard difference.
This would mean that over the course of the follow-up, an average of approximately 2 additional deaths will  occur for
each month of follow-up in each 100 persons under the control treatment alive 
at the start of the month and {\em who would also be alive under the IFN treatment}, compared with each 100 IFN treated
persons alive at the start of the month and  {\em who would also be alive under the control treatment}. As suggested before, the assumption that $T^0\cip T^1|V$ is implausible however.
%\footnote{Should we add the suggested sensitivity analysis for the HR as well?
%{\color{blue}TM: I have now added something on this, see above}
%}

\noindent

\subsection{Gastrointestinal tumour study}
%{\color{blue}
\cite{Stablein:1985aa}
presented survival data from a randomized 
clinical trial on locally unresectable gastric cancer. Half of the total 90 patients were assigned to chemotherapy, and the other half to combined chemotherapy and radiotherapy. 
%More details can be found in 
%\cite{GAST:1982aa}
 %Gastrointestinal Tumor Study Group (1982) "Gastrointestinal Tumor Study Group: Schein, P. D., Stablein, D. M., Bruckner, H. W., Douglass, H. O., Mayer, R., et al. (1982). A comparison of combination chemotherapy and combined modality therapy for locally advanced gastric carcinoma. Cancer 49, 1771-1777." 
It was suggested that there was 
 superior survival for patients who received  chemotherapy, but only in the first year or so. The same application was considered by
  \cite{collett:survival}  p. 386-389.
%There were two censored observations in the first treatment arm and six in the second.
For illustrative purposes we consider here the first 720 days of follow up corresponding to the two first time periods considered by 
  \cite{collett:survival}.
Figure 8 displays the Kaplan Meier curves corresponding to the two groups, and shows that the survival curves become close   at the end of the considered time interval. Applying a Cox regression model with  time-by-treatment-interaction, allowing for
  separate HR's before and after 1 year of follow-up gives 
 estimated HR's (Combined vs chemotherapy) of  
 2.40     (95\% CI; 1.25, 4.63) in the first year and  0.78   (95\% CI; 0.34,1.76) thereafter (the supremum score process test of Lin et al. (1993) gave no convincing evidence against these two models, with P 0.07 in the first interval and P 0.26  in the second interval).
Thus one might be tempted to conclude that 
 the chemotherapy is beneficial in the first year only, and that there might even be a reverse effect afterwards,
  see  \cite{collett:survival} for a similar analysis and conclusion.
 %Indeed one can argue .... see the  Gastrointestinal Tumor Study Group (1982) for  more details in this direction. 
 However, arguing based on the two hazard ratios is subtle as we have shown. If chemotherapy is more effective then there will be more and more frail subjects in that group making the interpretation of hazard function difficult. 
We illustrate this using the two estimated hazard ratios and taking the frailty variable to 
 be Gamma distributed mean and variance equal to
 % {\color{blue} 
 $\theta$  corresponding to a Kendall's $\tau$  of 0.3.
 %}
  Figure 9 displays the estimated $\lambda(t,A=1,Z)/\lambda(t,A=0,Z)$ which  is seen to 
 depend on time, being larger than
2.4 and increasing towards 3.7 in the first year and, after the change-point (1 year), starting at around 1.2 and then decreasing but being larger than 1 at all times.
% at all times in contrast to the naive conclusion based on the change point Cox model.
%}

\section{Concluding Remarks}

We have argued that the treatment effect in a proportional hazards model carries a causal interpretation, but that its interpretation is subtle. The proportional hazards assumption does not express, for instance, that treatment works equally effectively at all times, as the hazard ratio at a given time mixes differences between treatment arms due to treatment effect as well as selection. The danger of over interpreting hazard ratios become most pronounced when the hazard ratio is not constant over time (e.g. when the hazard ratio is below 1 for some time and then becomes 1). We have argued that this cannot be interpreted as implying that treatment effectiveness disappears after some time.
In our opinion, this is the source of much confusion, and a real concern. Non-constant hazard ratios are indeed fairly common in real life because the proportional hazards assumption is a rather unstable assumption in the following sense. Even when valid in some population, this assumption is likely to fail in subgroups of that population (e.g. if one studies men and women separately), and vice versa. This makes the assumption, at best, an approximation in practice.  

We have suggested possibilities to estimate hazard ratios that are causally interpretable because they compare intensities at a given time $t$ with and without treatment for the same patient population: those who would survive that time, no matter what treatment. Inferring such hazard ratios necessitates a sensitivity analysis, however. Furthermore, they have the disadvantage of describing the effect for an unknown subgroup of the population. A better strategy in practice, when interest lies in the dynamic aspects of a treatment, is therefore to design the study such that the collected data provide immediate insight into the dynamic aspects of treatment (e.g. by modifying treatment assignments over time).
 
 \bibliographystyle{apalike} % Tells it how you want references displaying in the bibliography.
\bibliography{CausalCox}

\begin{thebibliography}{}

\bibitem[Aalen et~al., 2015]{AalenCook15}
Aalen, O.~O., Cook, R.~J., and R{\o}ysland, K. (2015).
\newblock Does cox analysis of a randomized survival study yield a causal
  treatment effect?
\newblock {\em Lifetime Data Anal}, 21(4):579--93.

\bibitem[Collett, 2015]{collett:survival}
Collett, D. (2015).
\newblock {\em Modelling Survival Data in Medical Research}.
\newblock CRC Press, Boca Raton.

\bibitem[Hern{\'a}n, 2010]{Hernan2010}
Hern{\'a}n, M. (2010).
\newblock The hazards of hazard ratios.
\newblock {\em Epidemiology (Cambridge, Mass.)}, 21(1):13--15.

\bibitem[Lin et~al., 1993]{lin:wei:ying:1993}
Lin, D.~Y., Wei, L.~J., and Ying, Z. (1993).
\newblock Checking the {C}ox model with cumulative sums of martingale-based
  residuals.
\newblock {\em Biometrika}, 80:557--572.

\bibitem[Martinussen and Scheike, 2006]{ms06}
Martinussen, T. and Scheike, T.~H. (2006).
\newblock {\em Dynamic Regression Models for Survival Data}.
\newblock Springer-Verlag, New York.

\bibitem[Martinussen and Vansteelandt, 2013]{MartVan_Collap}
Martinussen, T. and Vansteelandt, S. (2013).
\newblock On collapsibility and confounding bias in cox and aalen regression
  models.
\newblock {\em Lifetime Data Anal}, 19(3):279--96.

\bibitem[Rubin, 2000]{10.2307/2669382}
Rubin, D.~B. (2000).
\newblock Causal inference without counterfactuals: Comment.
\newblock {\em Journal of the American Statistical Association},
  95(450):435--438.

\bibitem[Scheike et~al., 2015]{Scheike:2015aa}
Scheike, T.~H., Holst, K.~K., and Hjelmborg, J.~B. (2015).
\newblock Measuring early or late dependence for bivariate lifetimes of twins.
\newblock {\em Lifetime Data Anal}, 21(2):280--99.

\bibitem[Stablein and Koutrouvelis, 1985]{Stablein:1985aa}
Stablein, D.~M. and Koutrouvelis, I.~A. (1985).
\newblock A two-sample test sensitive to crossing hazards in uncensored and
  singly censored data.
\newblock {\em Biometrics}, 41(3):643--52.

\bibitem[Uno et~al., 2014]{Uno14}
Uno, H., Claggett, B., Tian, L., Inoue, E., Gallo, P., Miyata, T., Schrag, D.,
  Takeuchi, M., Uyama, Y., Zhao, L., Skali, H., Solomon, S., Jacobus, S.,
  Hughes, M., Packer, M., and Wei, L.-J. (2014).
\newblock Moving beyond the hazard ratio in quantifying the between-group
  difference in survival analysis.
\newblock {\em J Clin Oncol}, 32(22):2380--5.

\bibitem[Vansteelandt et~al., 2014]{Vansteelandt:2014aa}
Vansteelandt, S., Martinussen, T., and Tchetgen, E.~T. (2014).
\newblock On adjustment for auxiliary covariates in additive hazard models for
  the analysis of randomized experiments.
\newblock {\em Biometrika}, 101(1):237--244.

\bibitem[White and Royston, 2009]{WhiteRoyston_09}
White, I.~R. and Royston, P. (2009).
\newblock Imputing missing covariate values for the cox model.
\newblock {\em Stat Med}, 28(15):1982--98.

\bibitem[Zhao et~al., 2016]{Zhao_RMST}
Zhao, L., Claggett, B., Tian, L., Uno, H., Pfeffer, M.~A., Solomon, S.~D.,
  Trippa, L., and Wei, L.~J. (2016).
\newblock On the restricted mean survival time curve in survival analysis.
\newblock {\em Biometrics}, 72(1):215--21.

\end{thebibliography}

\nothere{
\noindent
%then $P(T>u | do(A=a,T>t))$ depends on $a$ (but that's sort of ok).
 A variation of this is to look at 
$$
P(T>u | do(A=a,T>t),A=\tilde a)
$$
%and see if this depends on $\tilde a$. It does not, regardless of which of the two DAG's we consider. 
When studying these two latter probabilities it is  easier to use the notation $S_t:=I(T>t)$ and then consider   DAG's 

 \begin{figure}[h!]
$$
\xymatrix{
 &S_t  \ar[r]&S_u\\
A  \ar[ur]\ar[urr]  & &
}
 $$
%  \caption{ Corresponds to $A\longrightarrow T$}
\end{figure}
\vspace{1cm}

\begin{figure}[h!]
$$
\xymatrix{
& Z\ar[rd]   &\\
A  \ar[rr] &  &T
}
 $$
%  \caption{}
\end{figure}

 \begin{figure}[h!]
$$
\xymatrix{
 Z  \ar[dr] \ar[drr]& &\\
 &S_t  \ar[r]&S_u\\
A  \ar[ur]\ar[urr]  & &
}
 $$
%  \caption{}
\end{figure}
}

\section*{Appendix A}
\subsection*{A.1 Binary frailty variable}

Let $Z$ be binary, for instance $P(Z=0.2)=0.2$ and $P(Z=1.2)=0.8$, low and high risk groups, then $E(Z)=1$, and  similar results as in the Gamma-distribution case (Section 3.1)
are obtained. In this  binary case, the Laplace transform is 
$$
\phi_Z(u)=0.2e^{-0.2u}+0.8e^{-1.2u}.
$$
Therefore
$$
\mbox{HR}_Z(t)\frac{g_Z(e^{-\Lambda(t,a=1)})}{g_Z(e^{-\Lambda(t,a=0)})}=\frac{\lambda(t;a=1)}{\lambda(t;a=0)},
$$
where
\begin{equation}
\label{g_Z}
g_Z(u)=\{D\log{(\phi_Z)}\}\{\phi_Z^{-1}(u)\},
\end{equation}
which is an increasing function. Hence, if we take $\beta_1<0$ and $\beta_2=0$, then again
$$
\mbox{HR}_Z(t)<1
$$
for all $t$.

\nothere{
Assume that the DGP is governed by a $\lambda(t;A,Z)$ with $A$ binary,  and $A$ independent of  $Z$. 
Further, assume that $T^0$ and $T^1$ are conditionally independent given $Z$.
Under suitable regularity conditions, it follows that
$$
\frac{\lambda^1(t)}{\lambda^0(t)}=\frac{E[\lambda(t;A=1,Z)\exp{\{-\Lambda(t; A=0,Z)-\Lambda(t; A=1,Z)\}}]}{E[\lambda(t;A=0,Z)\exp{\{-\Lambda(t; A=0,Z)-\Lambda(t; A=1,Z)\}}]}.
$$
Hence, 
\begin{equation}
\label{causal_hz}
\frac{\lambda^1(t)}{\lambda^0(t)}=\frac{\lambda(t;A=1,Z)}{\lambda(t;A=0,Z)}
\end{equation}
if 
$$
\lambda(t;A,Z)=g(t,A)h(t,Z)
$$
for some functions $g$ and $h$. If $\lambda(t;A,Z)=Z\lambda^*(t;A)$, as assumed in Section 2, then \eqref{causal_hz} holds.

We now suggest a way of estimating the $\mbox{HR}_Z(t)$ assuming we have knowledge about $Z$, for example matching subjects based on known risk factors, and randomize treatment within pairs. Let $(T_{0k},T_{1k})$ denote the $k$th pair with active treatment denoted by 1. Assume for the moment that there is no censoring but censoring due to end of study, and that subjects in a pair enter the study at the same point in time. The idea is that we censor $T_{ak}$, $a=0,1$, at $C_k=min(T_{0k},T_{1k})$. Let 
$N_{ak}(t)=I(T_{ak}\le t,T_{ak}\le C_k)$, and $Y_{ak}(t)=I(t\le C_k)$, $a=0,1$.  It can then be seen that
the intensity of $N_{ak}(t)$ is given by
$$
-Y_{ak}(t)\lambda^{\ast}(t;a)D\log{[\phi_Z\{\Lambda^{\ast}(t;a=0)+\Lambda^{\ast}(t;a=1)\}]},
$$
the ratio being, apart from the at risk indicators,
$$14
\frac{\lambda^{\ast}(t;a=1)}{\lambda^{\ast}(t;a=0)}=
\frac{Z_k\lambda^{\ast}(t;a=1)}{Z_k\lambda^{\ast}(t;a=0)}=HR_Z(t)
$$
so in this way we are targeting the true HR. I guess a stratified model (using info on the matched pairs), and allowing for time-changing treatment effect will also do the job. Both strategies suffer from the fact they rely on assuming the multiplicative structure.
}

\subsection*{A.2 Frailty model arising by marginalization}
It was shown in Section 3.1  that one can always pick a DGP $\lambda(t;A,Z)$ so that the Cox model holds marginally, only conditioning on the observed $A$.  Rename
$Z$ to $Z_1$. We show now that similarly we can also pick a DGP $\lambda(t;A,Z_1,Z_2)$ so that it marginalizes to  $\lambda(t;A,Z_1)$ that further  marginalizes to  $\lambda(t;A)$, the latter  being the Cox model. For ease of calculations, let 
$$
\lambda(t;A,Z_1,Z_2)=Z_2\lambda^*(t;A,Z_1)
$$
with $Z_2$ being Gamma distributed with mean and variance equal to 1, and independent of $Z_1$ and $A$. Similar calculations as those in Section 3.1 gives the following expression
\begin{equation}
\label{Z1Z2}
\lambda(t;A,Z_1,Z_2)=Z_1Z_2\lambda_0(t)e^{\beta A}\exp{\bigl [ \Lambda_0(t)e^{\beta A}+Z_1\{
\exp{\{\Lambda_0(t)e^{\beta A}}-1\}\bigr ]}
\end{equation}
Hence, if the DGP is governed by \eqref{Z1Z2} then model \eqref{Cox_AZ}, with $\nu=\infty$, 
and  the marginal Cox model, only conditioning on $A$, are also correctly specified.

\subsection*{A.3 Selection and Cox model}
Assume that the Cox model  $\lambda(t;A)$ is correctly specified.
%Assume that is perfectly described by the Cox model. 
Will there always be selection? The answer is yes.
%, unless the event time is determined without error, but then the Cox model cannot describe the data in the first place.
The Cox model induces randomness as
$$
\Lambda_0(T)=e^{-A\beta}V,
$$
where $V$ is exponentially distributed with mean 1. But then
\begin{align*}
E(V|T>t,A=a)
=E(V|V>e^{a\beta}\Lambda_0(t),A=a)
=1+e^{a\beta}\Lambda_0(t).
\end{align*}
If $e^{\beta}<1$ then 
$$
E(V|T>t,A=1)<E(V|T>t,A=0).
$$
%and similarly with the case, where $e^{\beta}>1$.$\hfill \square$
\bigskip

%{\color{blue}
\subsection*{A.4 A DGP with $\mbox{HR}(t)$ and   $\mbox{HR}_{Z}(t)$ being different}
Let   $Z$ be Gamma distributed with mean $1-\alpha$, $\alpha\in (0,1)$, and variance $1$, and let $V_0$ and $V_1$ be independent Gamma distributed with mean $\alpha$ and variance $1$.
Exposure $A$ is binary with $P(A=1)=1/2$. Generate data as follows:
$T^0=V_0+Z,$
$T^1=e^{-\beta}(V_1+Z)$ and
$T=(1-A)T^0+AT^1.$
Condition \eqref{cip} holds, and the marginal Cox model is also correctly specified.
 But in this case $\mbox{HR}(t)\neq \mbox{HR}_{Z}(t)$, and it also easily seen that $\mbox{HR}_{Z}(t)$ depends on $Z$.
  Different values of $\alpha$ results in different values of Kendall's $\tau$ thus controlling the correlation between $T^0$ and $T^1$.
  %Figure 2 shows   $\mbox{HR}(t)$ in a situation with Kendall's $\tau$ equal to 0.33
\nothere{
\begin{center}
\begin{figure}[h!]
  \includegraphics[width=14cm,height=12cm]{/Users/dkn606/Documents/Arbejde/papers/Cox_selection_bias?/Numerical/HR_t_ny1.pdf}
  %\caption{  True HR.}
     \caption{ {\small $\mbox{HR}(t)$ in scenarios with $e^{\beta}=0.5$ and with  Kendall's $\tau$ equal to 0.33.
     }}
 \end{figure}
  \end{center}
  }

%}

%{\color{blue}
\subsection*{A.5 Hazard differences} 
We assume \eqref{cip} and \eqref{ah}. Then 
$$
\lim_{h\rightarrow 0}P(t\le T^1< t+h|T^0\ge t,T^1\ge t)=\psi(t)+
\frac{E\{\omega(t,Z)e^{-2\int_0^t\omega(s,Z)\, ds}\}}{E\{e^{-2\int_0^t\omega(s,Z)\, ds}\}}
$$
and 
$$
\lim_{h\rightarrow 0}P(t\le T^0< t+h|T^0\ge t,T^1\ge t)=
\frac{E\{\omega(t,Z)e^{-2\int_0^t\omega(s,Z)\, ds}\}}{E\{e^{-2\int_0^t\omega(s,Z)\, ds}\}}
$$
and therefore 
$$
\psi(t)=\lim_{h\rightarrow 0}P(t\le T^1< t+h|T^0\ge t,T^1\ge t)-\lim_{h\rightarrow 0}P(t\le T^0< t+h|T^0\ge t,T^1\ge t).
$$
%}

\newpage

$$
\mbox{    }
$$

\begin{center}
\begin{figure}[h!]
  \includegraphics[width=12cm,height=7cm]{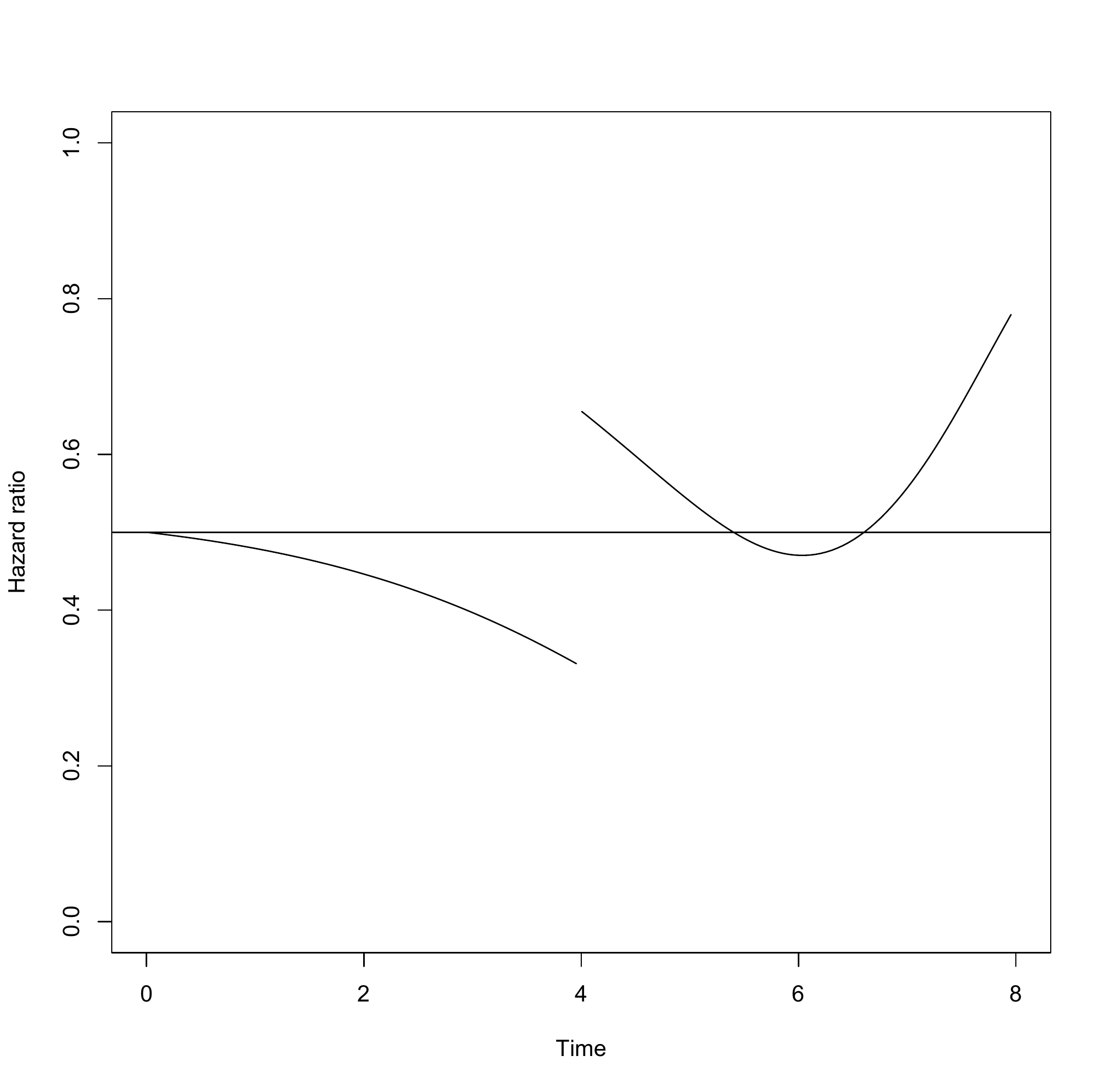}
  %\caption{  True HR.}
     \caption{ {\small Simulation study.  Plot of $\mbox{HR}_Z(t)$.}}
 \end{figure}
  \end{center}
  
\noindent
\begin{center}
\begin{figure}[h!]
  \includegraphics[width=12cm,height=7cm]{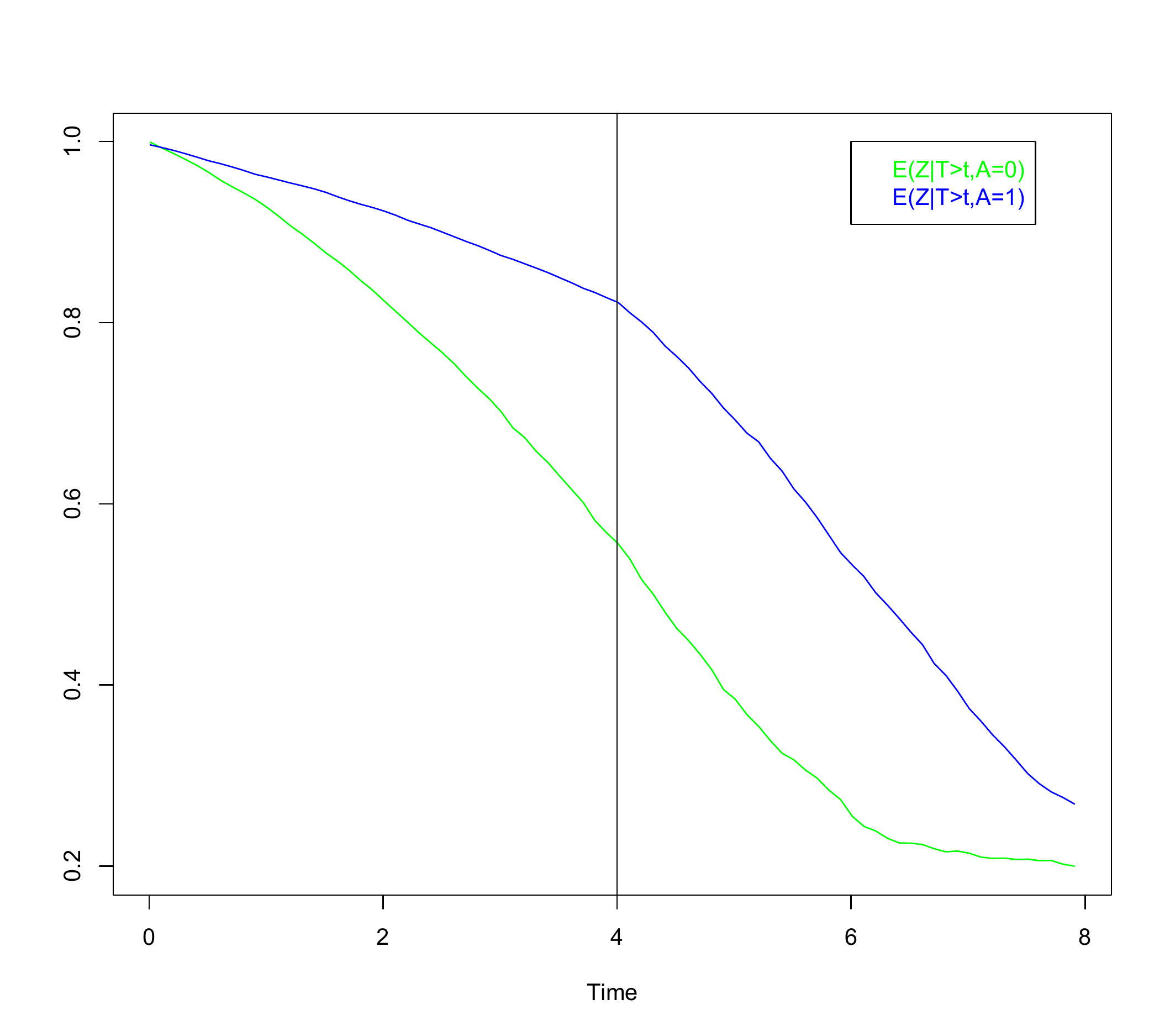}
  \caption{  {\small Simulation study.    Plot of $E(Z|T>t,A=a)$, $a=0,1$.}}
  \end{figure}
\end{center}

    \begin{center}
\begin{figure}[h!]
\label{Fig: HR_t}
  \includegraphics[width=14cm,height=12cm]{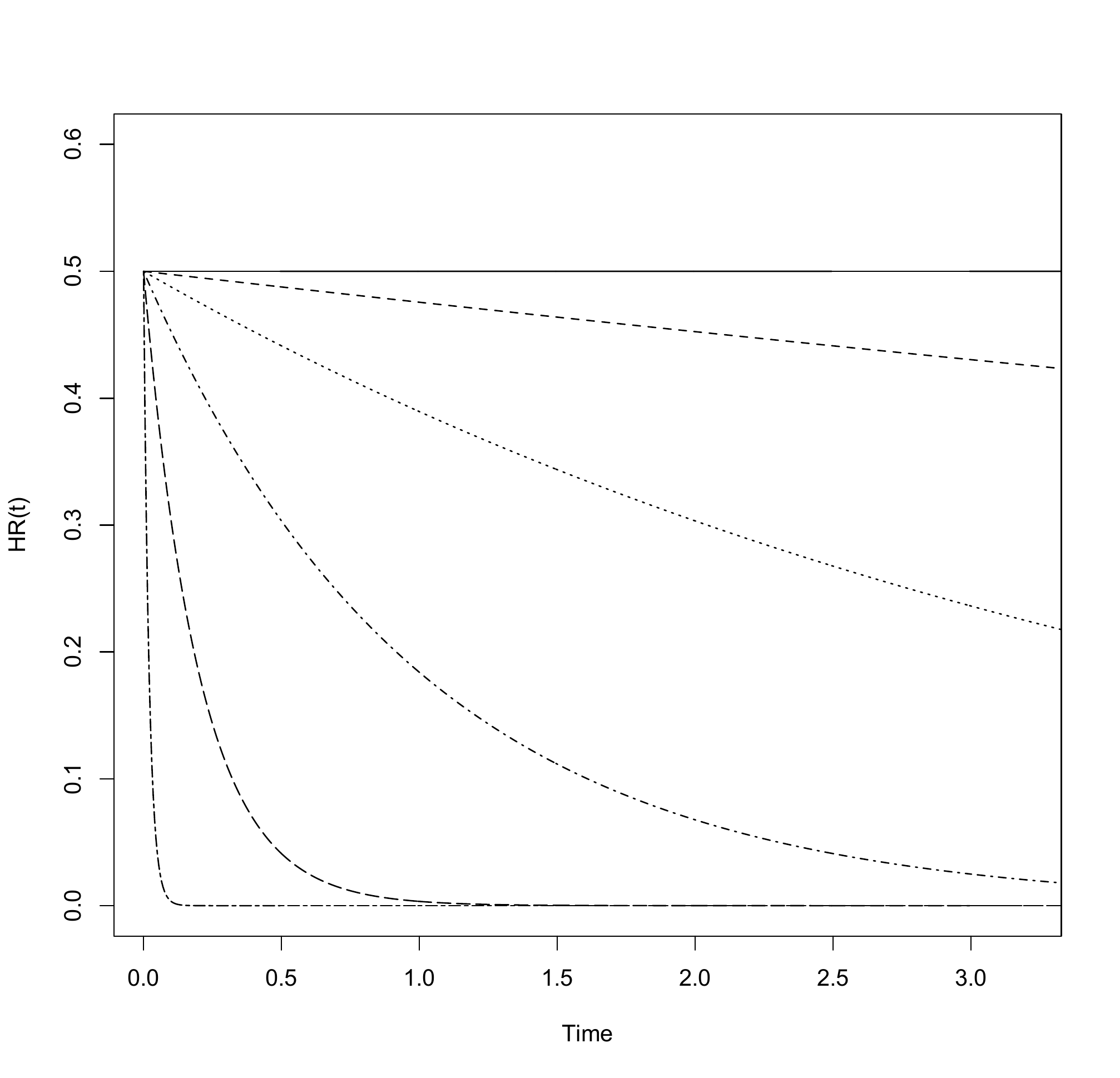}
  %\caption{  True HR.}
     \caption{ {\small $\mbox{HR}(t)$ in scenarios with $e^{\beta}=0.5$ and with  Kendall's $\tau$ equal to 0, 0.04, 0.2, 0.49, 0.83 and 0.98 (starting from top with $\tau$ equal  to 0, corresponding to independence between $T^0$ and $T^1$).
     }}
 \end{figure}
  \end{center}

\noindent
\begin{center}
\begin{figure}[h!]
  \includegraphics[width=12cm,height=7cm]{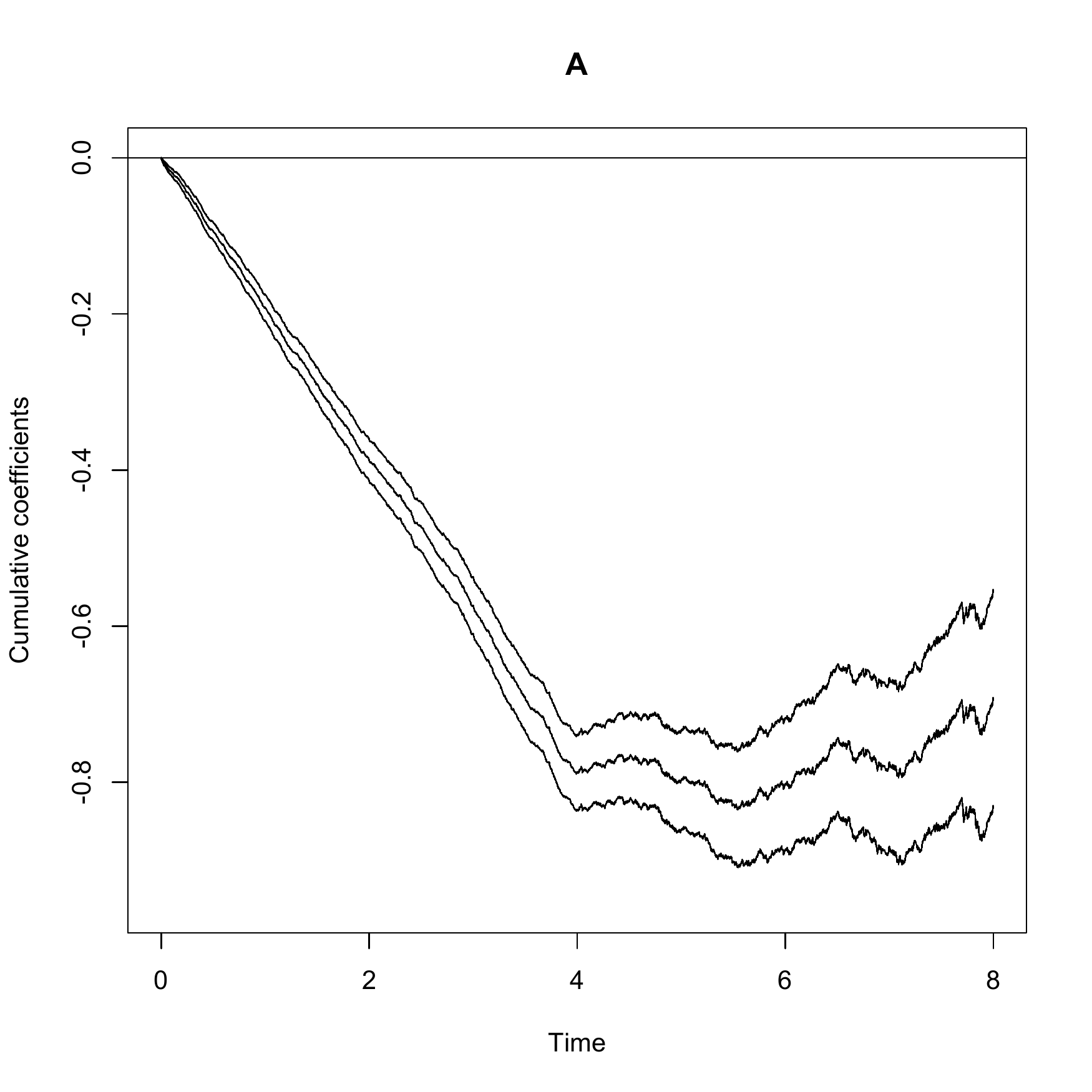}
  \caption{ {\small Simulation study.  Estimated cumulated regression coefficient and 95\% pointwise confidence bands obtained from fitting the Aalen additive hazards model.}}
  \end{figure}
\end{center}

%  \begin{center}
%\begin{figure}[h!]
 % \includegraphics[width=14cm,height=10cm]{HR_among_T0_T1_1.pdf}
%     \caption{ Simulation study.  Plot of  $\mbox{HR}(t)$.}
 %\end{figure}
  %\end{center}
  
  \begin{center}
\begin{figure}[h!]
\label{Fig: KM}
  \includegraphics[width=12cm,height=7cm]{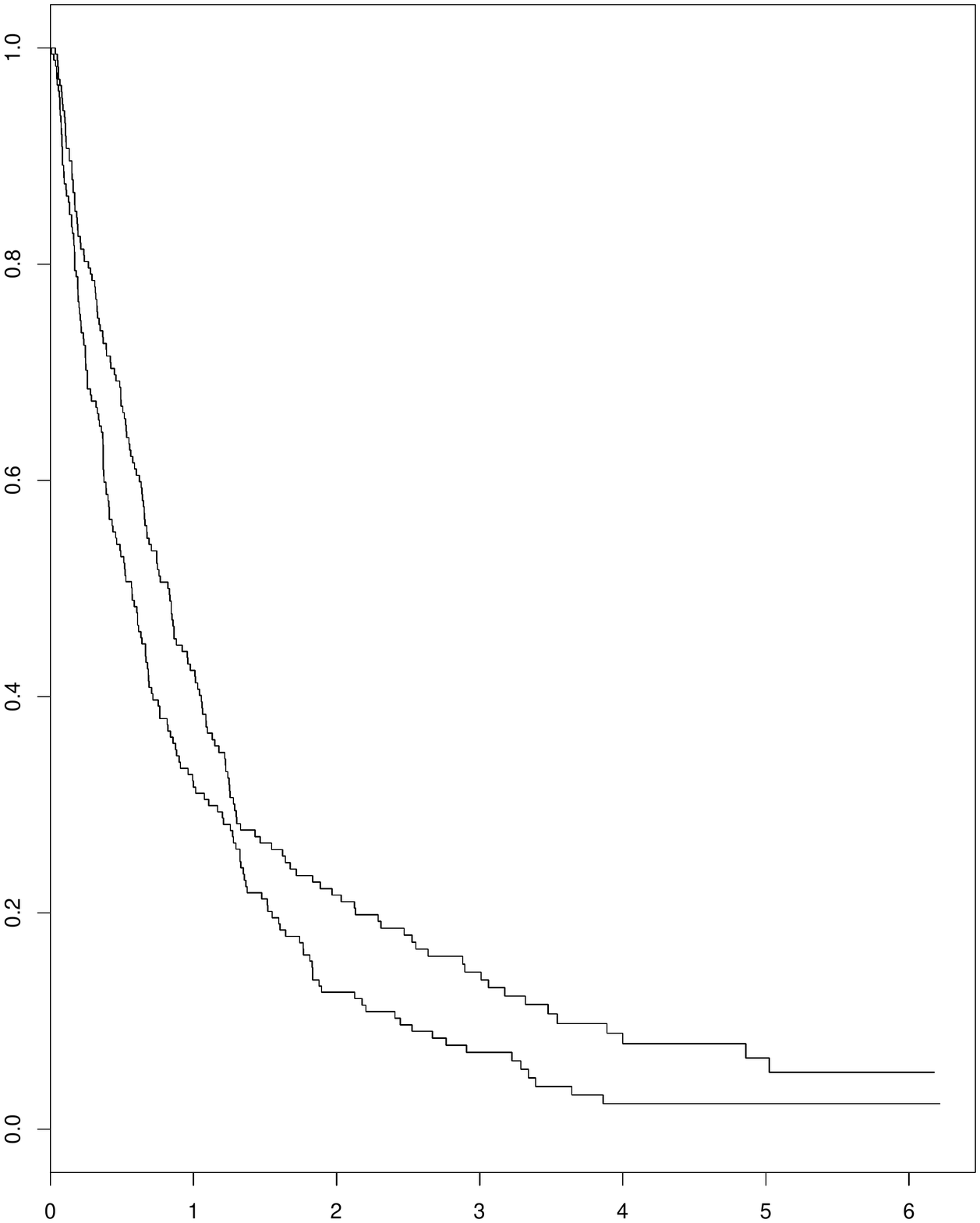}
  %\caption{  True HR.}
     \caption{ {\small MRC RE01 study. Kaplan Meier plot, control group (green curve) and IFN group (blue curve).}}
 \end{figure}
  \end{center}

   \begin{center}
\begin{figure}[h!]
\label{Fig: HR.t}
  \includegraphics[width=12cm,height=7cm]{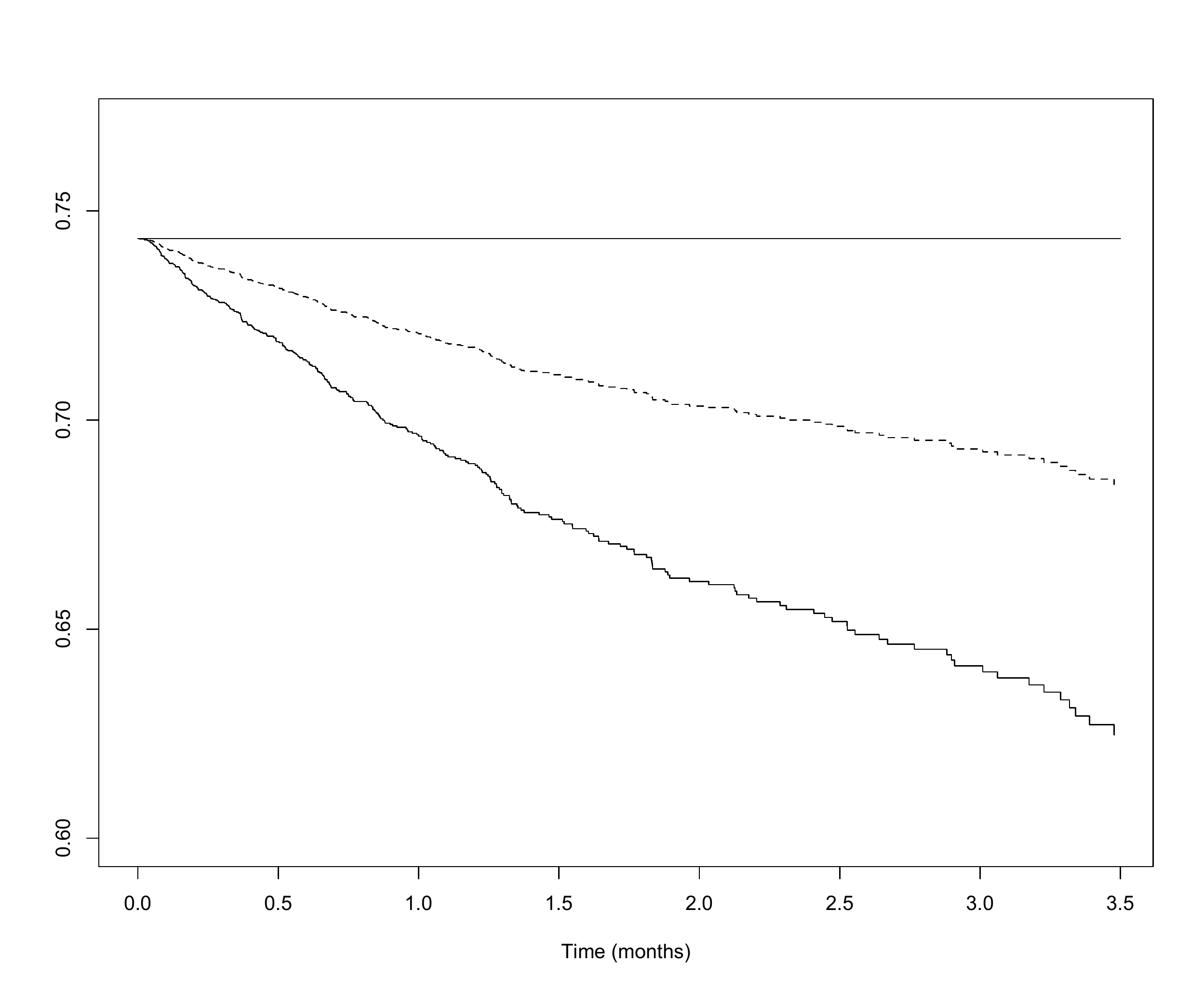}
  %\caption{  True HR.}
     \caption{ {\small MRC RE01 study. Estimated $\mbox{HR}(t)$ with Kendall's $\tau$ equal to 0.3 (dotted curve), 0.2 (broken curve) and 0.1 (full curve).
     Horizontal line corresponds to the Cox hazard ratio of 0.74.
     }}
 \end{figure}
  \end{center}

  \begin{center}
\begin{figure}[t!]
\label{Fig: RR}
  \includegraphics[width=12cm,height=7cm]{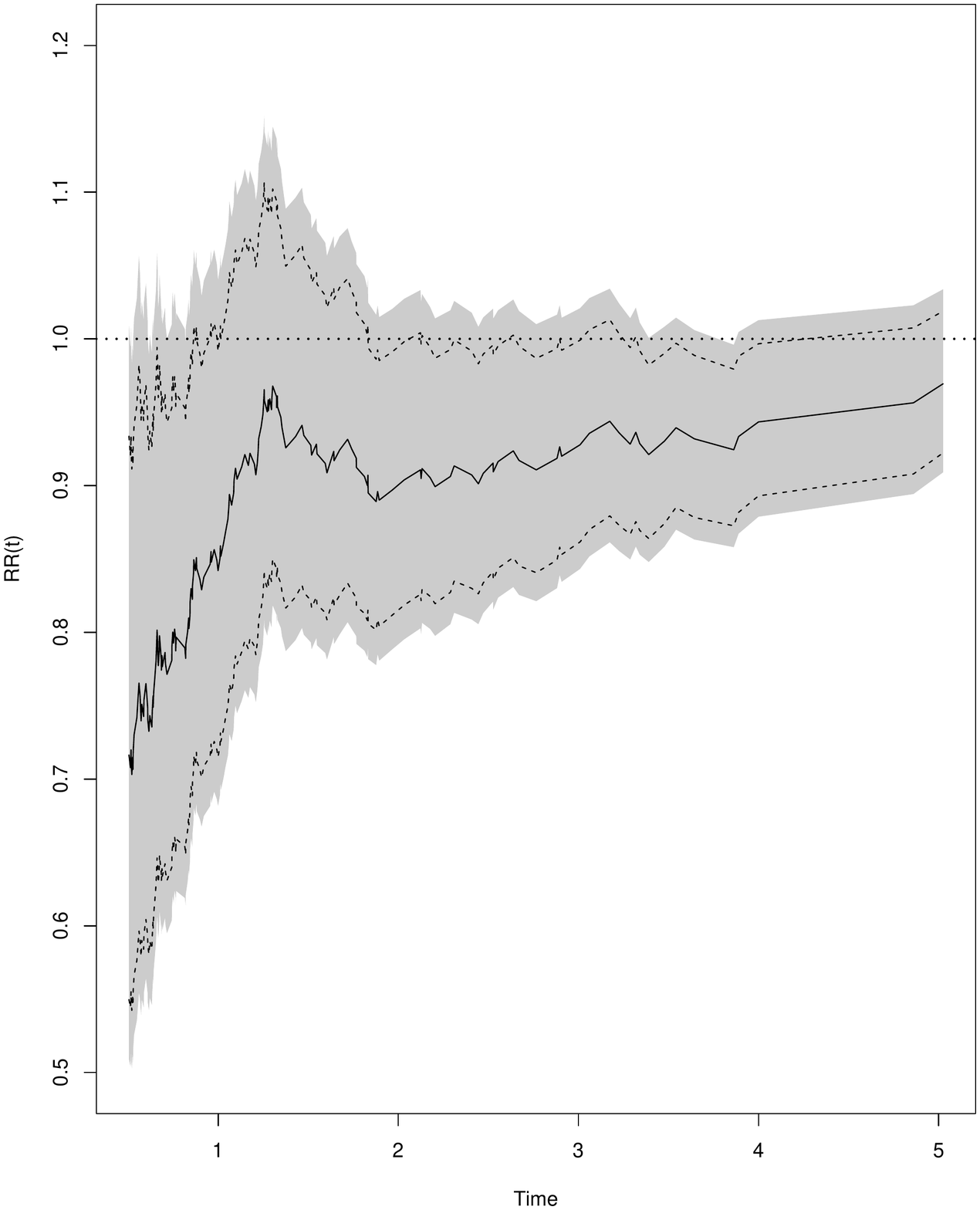}
  %\caption{  True HR.}
     \caption{ {\small MRC RE01 study. IFN treatment vs control treatment. Estimate of relative risk $\mbox{RR}(t)$ along with 95\% pointwise confidence bands (dashed curves) and 95\% uniform bands (shaded area).}}
 \end{figure}
  \end{center}

  %%% Aalen plots  not in 

 \begin{center}
\begin{figure}[h!]
\label{Fig: Gas_KM}
  \includegraphics[width=12cm,height=7cm]{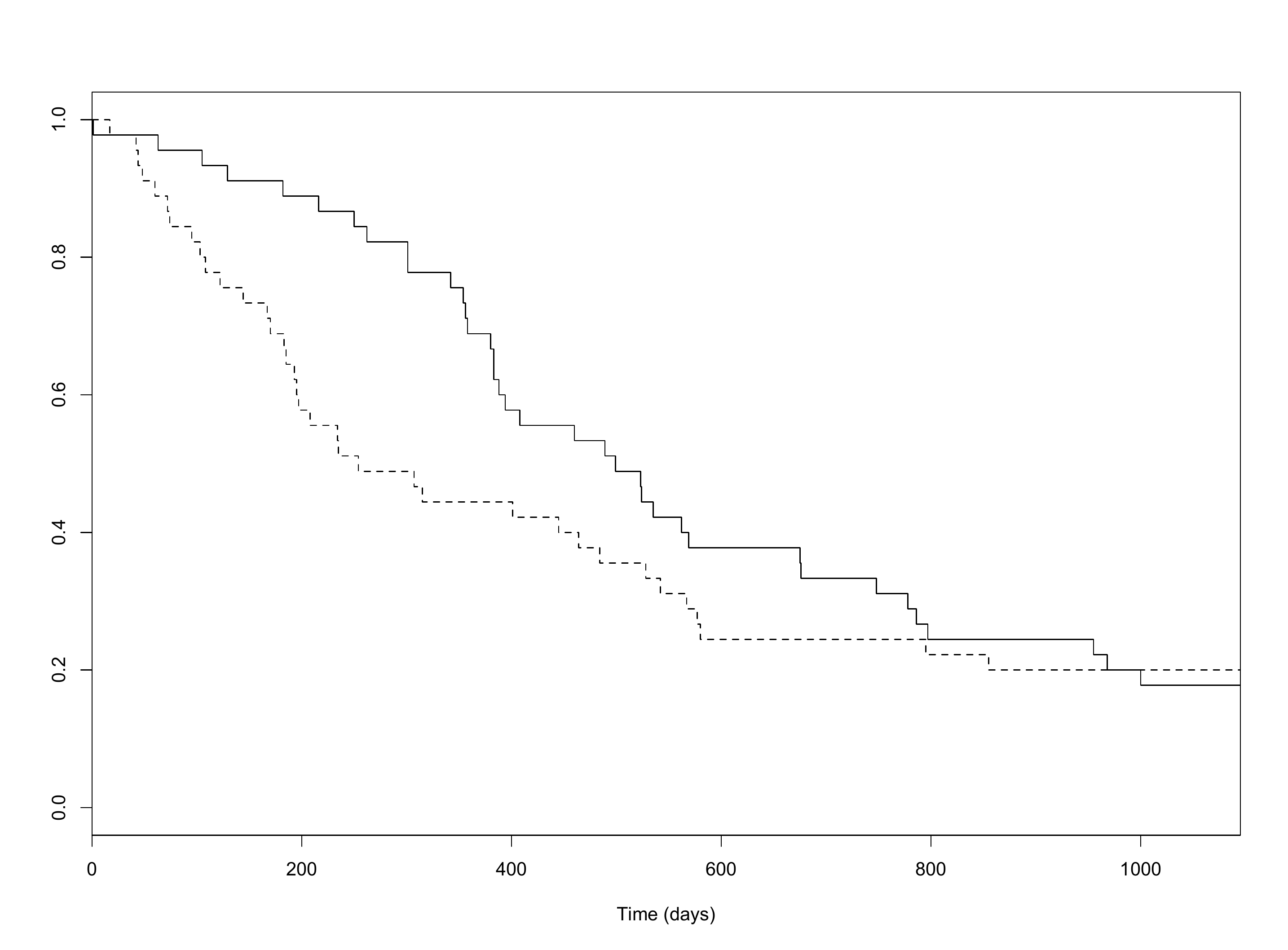}
  %\caption{  True HR.}
     \caption{ {\small Gastrointestinal tumour study. Kaplan-Meier plot, chemotherapy (full curve) and combined therapy (broken curve).
     }}
 \end{figure}
  \end{center}

 \begin{center}
\begin{figure}[h!]
\label{Fig: Gas_HR}
  \includegraphics[width=12cm,height=7cm]{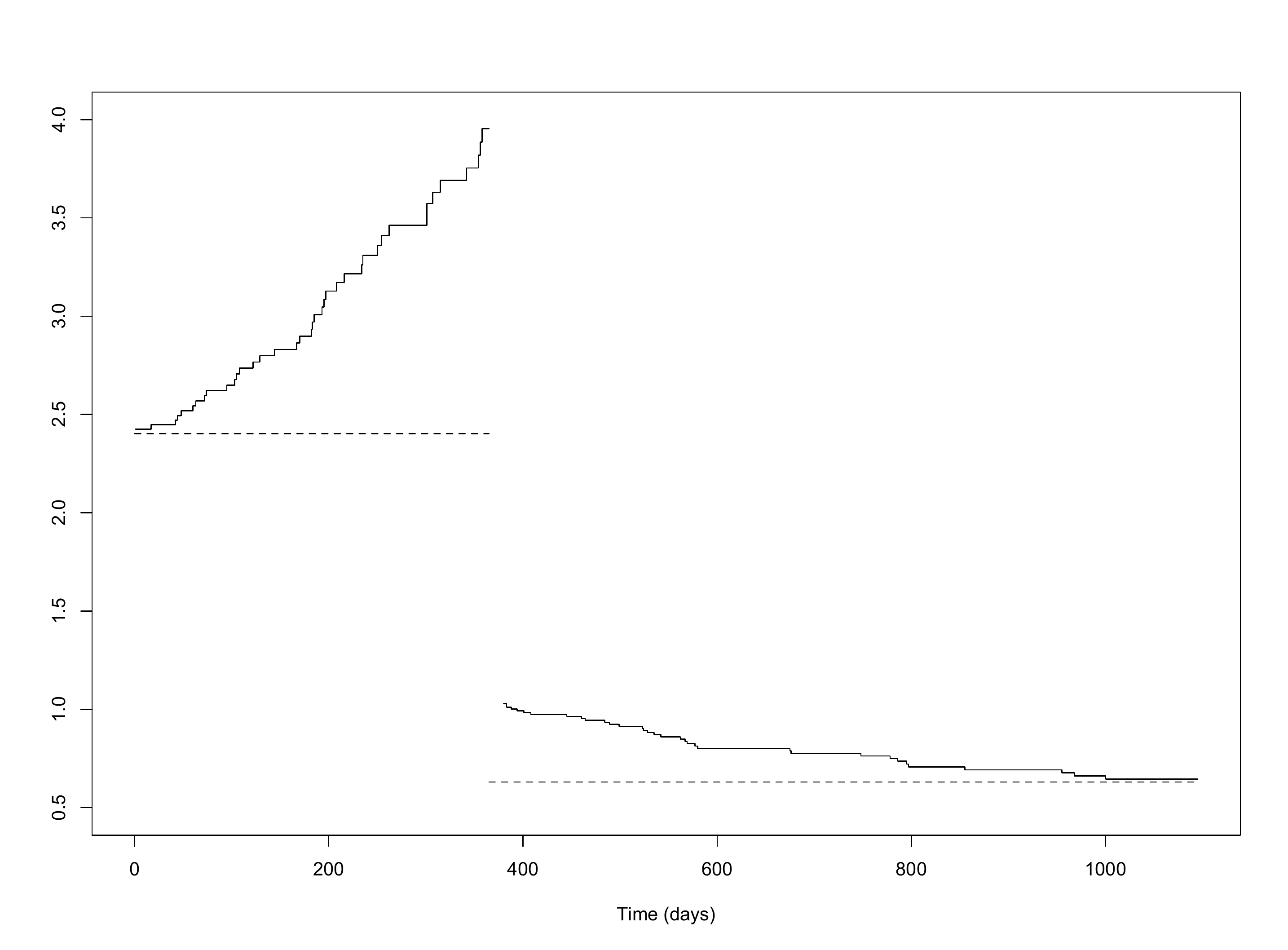}
  %\caption{  True HR.}
     \caption{ {\small Gastrointestinal tumour study. Plot of $\lambda(t;A=1,Z)/\lambda(t;A=0,Z)$ based on $\eqref{Cox_AZ}$ with change-point at 1 year and the frailty variable being Gamma distributed with mean 1 and variance so that Kendall's $\tau$ is equal to 0.3. Dashed curves show the estimated regression coefficients based on the change-point Cox analysis.
     Combined therapy vs chemotherapy.
     }}
 \end{figure}
  \end{center}

\end{document}